\documentclass[10pt,a4,twoside,hidelinks,rm]{article}
\usepackage[normalem]{ulem}
\usepackage{cancel}
\usepackage[vcentermath]{youngtab}
\usepackage{young}
\usepackage{ytableau}
\usepackage{theorem}
\usepackage{amssymb}
\usepackage[widespace]{fourier}
\usepackage{fancyhdr}
\usepackage{float}
\usepackage[bb=esstix,bbscaled=.5,cal=cm,calscaled=.92]{mathalfa}
\usepackage{bbm}
\usepackage{graphicx}
\usepackage[ansinew,latin1]{inputenc}
\usepackage[active]{srcltx}
\usepackage{amsfonts}
\usepackage{nccmath}
\usepackage{titlesec}
\usepackage{pb-diagram}
\usepackage[colorlinks=false]{hyperref}
\usepackage{color}
\usepackage{epsfig}
\usepackage[sans]{dsfont}
\usepackage{multicol}
\usepackage[left=4.5cm,top=2cm,right=4.5cm,bottom=2cm]{geometry}
\usepackage{amsmath}
\usepackage{yfonts}
\usepackage{pb-diagram}
\usepackage[refpage]{nomencl}
\usepackage{ifthen}
\usepackage{lipsum}
\usepackage{commath}

\renewcommand\thesection{\arabic{section}} 
\renewcommand\thesubsection{\arabic{section}.\arabic{subsection}} 
\titleformat{\section}[block]{\scshape\centering}{\thesection.}{1em}{} 
\titleformat{\subsection}{\bfseries}{\thesubsection.}{1em}{} 

\newcommand{\seqN}{{\rm seq}_N}
\newcommand{\seqr}{{\rm seq}_r}

\newcommand{\kk}{\mathcal{K}}

\newcommand{\C}{\mathbb C}

\newcommand{\BB}{\mathbb{B}}
\newcommand{\Comp}{{\mathcal Comp}_n}
\newcommand{\PP}{{\mathcal P}}
\newcommand{\Par}{{\mathcal Par}_n}

\newcommand{\Y}{\mathcal{Y}}
\newcommand{\YY}{\mathcal{Y}_{r,n}(q)}

\newcommand{\End}{{\rm End}}

\newcommand\es{\mathbbm{s}}
\newcommand\et{\mathbbm{t}}

\newcommand\diez{{10}}
\newcommand\once{{11}}
\newcommand\doce{{12}}
\newcommand\trece{{13}}
\newcommand\bs{\mathbf{s}}

\newcommand\bu{\mathbf{u}}
\newcommand\bv{\mathbf{v}}

\newcommand{\s}{\mathfrak{s}}
\newcommand{\U}{\mathfrak{u}}
\newcommand{\V}{\mathfrak{v}}
\newcommand{\T}{\mathfrak{t}}

\newcommand{\bT}{\pmb{\mathfrak{t}}}
\newcommand{\Bs}{\pmb{\mathfrak{s}}}

\newcommand{\Bv}{\pmb{\mathfrak{v}}}

\newcommand{\MP}{{ {Par}}_{r,n}}

\newcommand{\MC}{{ {Comp}}_{r,n}}
\newcommand{\MCN}{{ {Comp}}_{r,n, \le N}}
\newcommand{\MPN}{{ {Par}}_{r,n, \le N}}

\newcommand{\Si}{\mathfrak{S}}
\newcommand{\std}{{\rm Std}}
\newcommand{\rstd}{{\rm RStd}}
\newcommand{\Tab}{{\rm Tab}}

\newcommand{\PTLthree}{{\cal PTL}_3(q)}
\newcommand{\PTLfour}{{\cal PTL}_4(q)}
\newcommand{\PTLfive}{{\cal PTL}_5(q)}
\newcommand{\PTL}{{\cal PTL}_n(q)}
\newcommand{\PTLa}{{\cal PTL}^{\alpha}_n(q)}

\newcommand{\HH}{ \mathcal{H}_n(q)}

\newcommand{\E}{ {\mathcal E}_n(q)}
\newcommand{\Ea}{ {\mathcal E}_n^{\alpha}(q)}

\newcommand{\ETL}{ {\cal ETL}_{n,N}(q)}

\newcommand{\EaTLK}{ {\cal ETL}_{n,N}^{ {\cal K}, \alpha}(q)}
\newcommand{\EaTL}{ {\cal ETL}_{n,N}^{  \alpha}(q)}
\newcommand{\EaTLtwo}{ {\cal ETL}_{n,2}^{  \alpha}(q)}

\newcommand{\Eak}{ {\mathcal E}_n^{ {\cal K}, \alpha}(q)}

\newcommand\bS{\Sigma}
\newcommand\blambda{{\boldsymbol\lambda}}

\newcommand\be{\mathbb{E}}
\newcommand\bmu{{\boldsymbol\mu}}

\theorembodyfont{ } \theoremstyle{marginbreak}

\theoremstyle{plain}

\newtheorem{teo}{Theorem}
\newtheorem{coro}[teo]{Corollary}

\newtheorem{defi}[teo]{Definition}

\newtheorem{lem}[teo]{Lemma}
\newtheorem{propos}[teo]{Proposition}

\newenvironment{demo}
{\textsc{Proof.}} {\quad \hfill $\Box$}

\begin{document}

\title{\bf On the annihilator ideal in the $bt$-algebra of tensor space  }

\author{Steen Ryom-Hansen\thanks{Supported in part by FONDECYT grant 1171379  } }
\date{}   \maketitle

\begin{abstract}
\noindent \textsc{Abstract. }
We study the representation theory of the \textit{braids and ties} algebra,
or the $bt$-algebra, 
$ \E$.
Using the cellular basis $ \{ m_{\es\et} \} $ for $ \E$ obtained in previous joint work with
J. Espinoza we introduce two kinds of
permutation modules $ M(\blambda ) $ and $ M(\Lambda ) $ for $ \E$.
We show that the tensor product module $ V^{\otimes n } $ for $ \E $ is a direct sum of $ M(\blambda)$'s.
We introduce the dual cellular basis $ \{ n_{\es\et} \} $ for $ \E $ and study its action on
$ M(\blambda ) $ and $ M(\Lambda ) $.
We show that the annihilator ideal $ \cal I $ in $ \E $ of $ V^{\otimes n } $
enjoys a nice compatibility property with respect to 
$ \{ n_{\es\et} \} $.
We finally study the quotient algebra $ \E/\cal I $, showing in particular that it is a simultaneous generalization of
H\"arterich's 'generalized Temperley-Lieb
algebra' and Juyumaya's 'partition Temperley-Lieb algebra'.
\end{abstract}


\medskip
\noindent
MCS2010: 33D80.

\section{Introduction}

\medskip
In this paper we study the representation theory
of the \textit{braids and ties algebra} $ { \cal E}_{n}(q) \, $, or the \textit{bt-algebra} for short.
It was introduced by Aicardi and Juyumaya in \cite{AJ1}, via an \textit{abstraction} of the presentation 
for the Yokonuma-Hecke algebra $ \YY$ of type $ A_{n-1}$. The $bt$-algebra carries a diagrammatics
consisting of {\it braids and ties}, as illustrated below where braids are black
and ties are dashed red.
\begin{equation}\label{bt}
  \raisebox{-.45\height}{\includegraphics[scale=0.4]{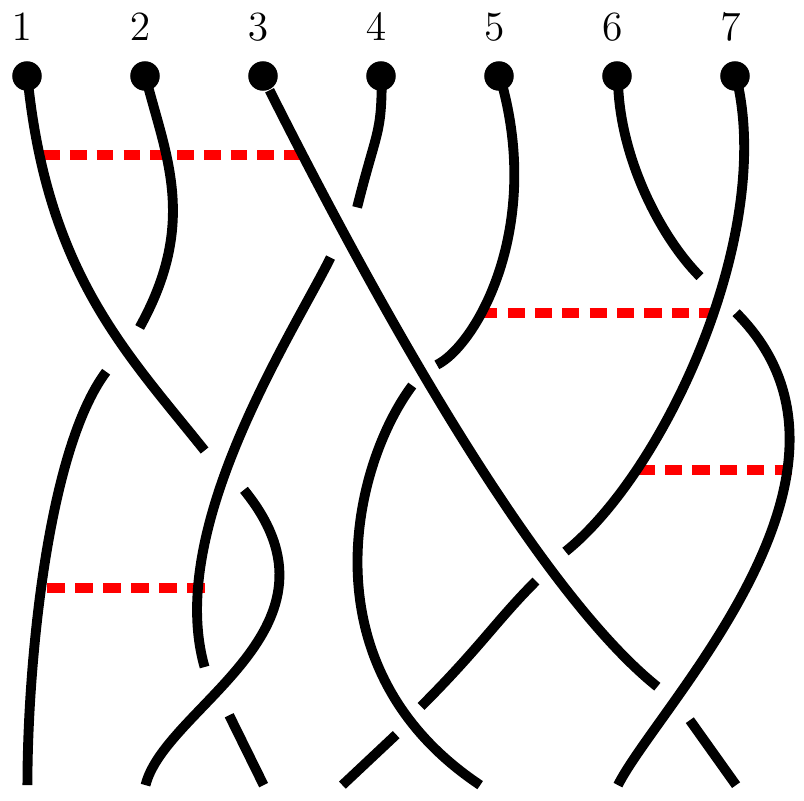}}
\end{equation}

The study of $\E $ has attracted much attention recently, 
both from knot theorists and from representation theorists,
see for example \cite{ArJu}, \cite{AJ2}, \cite{AJ3}, \cite{AJ4}, \cite{Ba},
\cite{ChlouPou},
\cite{ER},\cite{Flo}, \cite{JacondAn}, \cite{Juyumaya}, \cite{Marin},  \cite{Ry} and the references therein.
Interesting enough,
even though $\E $ is derived from $\YY$, both the knot theory and  
the representation theory of $\E$ are 
almost completely independent from those of $\YY$, in the sense 
that almost no statement in the literature about $\E $
is obtained directly from a corresponding statement
about $ \YY$.

\medskip
In our paper \cite{Ry}
we constructed a tensor space module $ V^{\otimes n} $ for
$\E $, extending Jimbo's classical tensor space module for the Iwahori-Hecke algebra
$ \HH $ of type
$ A_{n-1} $, 
and showed that it is faithful when the dimension of $ V $ is large enough, although not in general. Using this we
obtained a basis for $ \E $ on the form $ \{E_{A} g_w \}$ where
$A $ runs over 
{\it set partitions} on 
$ \mathbf{n}:=  \{ 1,2,\ldots, n \} $ and $ w $ over the symmetric group $ \Si_n$. 
In particular the dimension of $ \E $ is $b_n n! $, where $ b_n $ is the $n$'th Bell
number, that is the number of set partitions on $ \mathbf{n}$.

\medskip
The $ \{E_{A} g_w \}$-basis plays an important role in the applications of $ \E $ to knot theory as 
a key ingredient in the construction of a Markov trace on $ \E$, see for example
\cite{AJ2}, \cite{AJ3} and \cite{AJ4}. {\color{black} On knots, the resulting
invariant is equivalent to  
the HOMPLYPT-polynomial, but on links it is stronger than the HOMPLYPT-polynomial, see the discussion in
the Addendum to \cite{AJ2}}.

\medskip
In the present paper we develop methods that allow us to treat the $ \E $-structure on $ V^{\otimes n} $
when $ V^{\otimes n} $ is not a faithful, and in particular to determine 
the annihilator ideal $ \cal I $ for the $ \E $-action on $ V^{\otimes n} $
in all cases. 
Our methods rely on \textit{cellular algebra} theory, in particular on the cellular basis $ \{m_{\es \et} \} $
for $ \E $ that was introduced in \cite{ER}, in joint work with J. Espinoza, and hence they are
independent of the ground ring. 
We acknowledge that the construction of the $ \{m_{\es \et} \} $-basis for $ \E$ is combinatorially
much more involved than the construction of the $ \{E_{A} g_w \}$-basis, or other similar
non-cellular bases, 
but throughout 
our arguments rely crucially on the compatibility between  
the multiplicative structure and the order relation on $ \E $, given by the cell datum. In fact,
it appears to be very difficult to obtain the results of our paper without using an appropriate order
relation, even for the ground field $ \mathbb C $ and generic $ q $.

\medskip
The basis $ \{m_{\es \et} \} $ is a generalization of Murphy's 
\textit{standard basis} $ \{ x_{\s \T}\}$ for $ \HH$. 
In the classical representation theory of $ \HH $, see for example \cite{DJ} and \cite{DJM}, 
an important application of $\{ x_{\s\T}\}$ is to
introduce \textit{permutation modules} $ M(\lambda) $ for $ \HH$, where $ \lambda $ runs over 
\textit{integer partitions} of $ n $, and to realize Jimbo's tensor space module $ V^{\otimes n} $
as a direct sum of these $ M(\lambda)$'s. The $ M(\lambda)$'s
have many important properties, but for us it is of special relevance 
that they are endowed with canonical symmetric
bilinear forms $ (\cdot, \cdot)_{\lambda}$, that
satisfy a certain non-vanishing property with
respect to the \textit{dual basis} $ \{ y_{\s \T}\}$ of $ \{ x_{\s \T}\}$.

\medskip
In the paper we generalize
the $ M(\lambda) $'s to $ \E$-permutation modules $ M(\blambda) $, for $ \blambda $
running over \textit{multipartitions}, 
and show in our Theorem \ref{17}
that the $ \E$-tensor module $ V^{\otimes n} $ from \cite{Ry} is 
a direct sum of these $ M(\blambda) $'s in analogy with the $ \HH$-case. 
The $ M(\blambda) $'s are endowed with bilinear forms, in analogy with the $ M(\lambda)$'s, 
but the parametrizing poset $ {\cal L}_n $ for $ \E $ is a combinatorial object
which is more complicated than 
multipartitions, 
and so the $ M(\blambda) $'s are not sufficient for our determination of
the annihilator ideal
$ \cal I$.

\medskip
To solve this problem we introduce, via $ \{ m_{\es \et}\}  $, another kind of permutation modules
that we denote
$ M(\Lambda)$. 
We show
in a series of Lemmas, culminating with Lemma
\ref{samespirit}, 
that the canonical 
bilinear form
$ (\cdot, \cdot)_{\Lambda}$ on $ M(\Lambda)$
indeed satisfies the nonvanishing property
with respect to the dual basis $ \{ n_{\es \et}\}$.

\medskip
The above results constitute the bulk of work of our paper.
The situation is not completely analogous to the situation 
for $ \HH$, 
but in general there is an inclusion $ M(\Lambda) \subseteq M(\blambda) $ 
and this is in fact enough 
to execute 
H\"arterich's wonderful compatibility argument, see \cite{Har}, for the determination of the annihilator ideal of the action
of $ \HH $ in Jimbo's tensor space $ V^{\otimes n} $, but in the $\E$-setting. We do so in Theorem \ref{mainT}. 

\medskip
When working over a ground field $ \kk $, 
Theorem \ref{mainT} gives a $ \kk$-vector space basis for $ \cal I $ in terms of an explicit
subset of $ \{ n_{\es \et }\} $, thus generalizing H\"arterich's compatibility result. 
Note that
H\"arterich's 
compatibility result is the starting point for the author's joint work with D. Plaza on KLR-gradings
on the Temperley-Lieb algebras of type $A $ and $ B$, see \cite{PlaRy}.
Note also that other compatibility results for $ \{ m_{\es \et } \} $ were obtained
in \cite{ER}.

\medskip
In the final paragraphs of the paper we introduce and study the quotient algebra $ \ETL: =\E/\mathcal I $.
It is a simultaneous generalization of H\"arterich's 'generalized Temperley-Lieb' algebra
and of Juyumaya's 'partition Temperley-Lieb algebra $\PTL$', that themselves are generalizations of the original
Temperley-Lieb algebra. Our Theorem \ref{mainT} allows us to determine the dimension of $ \ETL$
as the cardinality of an explicit subset of $ \{ n_{\es \et} \} $. This shows in particular that Juyumaya's
conjectural basis for $ \PTL $ is wrong, but gives a correct basis that may be useful for constructing a
Markov trace on $\PTL$. 

\medskip
Let us now indicate the layout of the paper. In the next section we introduce the basic notation
that shall be used throughout the paper. It is mostly concerned with combinatorial notions
related to the symmetric group $ \Si_n $, that is partitions, Young diagrams, Young tableaux, etc. 
Throughout the paper, $ \Si_n $ is viewed as a Coxeter group and we only consider 
subgroups of $ \Si_n $ that 
are  
parabolic subgroups, in the sense of Coxeter groups. It is well known that parabolic subgroups
of Coxeter groups give rise to distinguished coset representatives. Throughout we need these
distinguished coset representatives
and their associated decompositions of the group elements, 
and we therefore briefly explain 
why they exist, but from an algorithmic point of view which is possibly less known.
In the final paragraphs of this section we recall the definition of $ \E$, together with its $ \{ E_A g_w \} $ basis.

In section 3 we first recall from \cite{ER} the ideal decomposition 
\begin{equation}\label{deideal} \E=\bigoplus_{\alpha\in\Par} \Ea \end{equation} that
reduces the study of $ \E $ to the study of its summands $ \Ea$. 
We next recall from \cite{ER} the construction of the cellular basis $ \{ m_{\es \et} \} $
for $ \Ea $ and at the same time we introduce the new dual cellular basis $ \{ n_{\es \et} \}$ 
by making the appropriate adaptions.
It is a generalization of Murphy's dual standard basis 
$ \{ y_{\s \T} \} $ for $ \HH$. 
The cell poset $ {\cal L}_n(\alpha) $ for both bases consists of pairs $ (\blambda \mid \bmu) $
of multipartitions, satisfying certain conditions.

In section 4 we first introduce
the permutation module $ M(\blambda) $.
In Lemma \ref{descr} we describe a basis $\{ x_{\Bs } \} $ for $ M(\blambda) $ together with
its $ \Ea $-action.  We next introduce, for $ \Lambda \in  {\cal L}_n(\alpha) $, the permutation module
$ M(\Lambda) $. We describe in the Lemmas \ref{descr2} and
\ref{descr2A} 
a basis $\{ m_{\Bs } \} $ together with the
$ \Ea $-action on it. We then introduce the bilinear form $ ( \cdot, \cdot)_{\Lambda} $ on
$ M(\Lambda) $ and show that it is symmetric and invariant. Finally, in Lemma
\ref{samespirit}
we prove the following crucial property, alluded to above
\begin{equation} ( m_{ \et } n_{ \et^{\prime} \es^{\prime}     } , m_{\es})_{\Lambda} \neq 0.
\end{equation}  

In section 5 we obtain the main results of our paper. We first give the decomposition
of $ V^{\otimes n} $ corresponding to (\ref{deideal}) and then the decomposition
in terms of the $ M(\blambda) $'s.
Finally, in the main 
Theorem \ref{mainT}, we give the description of $ \cal I$, proving that
a subset of $ \{ n_{\es \et} \} $ induces a basis for it.

\medskip
It is a great pleasure to thank J. Espinoza for many useful discussions related to the contents of the
paper, for sending us comments on preliminary versions of it, and for sharing with us his
formula for $ \dim \PTL $.
We also
thank J. Juyumaya for communicating to us his and P. Papi's MAGMA calculations
and for pointing out an error in a preliminary version of the paper.
{\color{black}We thank the referee for useful comments that helped us improve the text.}
Finally, it is a special pleasure to 
thank M. H\"arterich for introducing us to Murphy's standard basis
and all its deep properties.

\section{Notation and Basic Concepts}
We fix the ground ring $S:=\mathbb{Z}[q,q^{-1}]$,
where
$q$ is an indeterminate.


\medskip
Let $\Si_n$ be the symmetric group consisting of bijections of
$ \mathbf{n}:=  \{ 1,2,\ldots, n \} $.
Throughout we consider the natural action of $ \Si_n $ on $ \mathbf{n} $ as a right action.
As is well known, $\Si_n$ is a Coxeter group on $\bS_n:=\{s_1,\ldots, s_{n-1}\} $ with relations
\begin{alignat}{3}
s_is_j&=s_js_i &&\quad\mbox{ for }\; |i-j|>1\label{s1}\\
s_is_{i+1}s_i&=s_{i+1}s_is_{i+1}&&\quad\mbox{ for } i=1,2,\ldots,n-2 \label{ss11}\\
s_i^2&=1&&\quad\mbox{ for } i=1,2,\ldots,n-1\label{sss111}
\end{alignat}
where $ s_i := (i,i+1) $.
We denote by $ < $ the Bruhat-Chevalley order 
and by $\ell(\cdot) $ the length function on $ \Si_n$.

\medskip
Let $\mathbb{N} $ denote the set of natural numbers and set $ \mathbb{N}^0 := \mathbb{N} \cup \{ 0 \}  $.
{\color{black}{For the notions of \textit{compositions} and \textit{partitions} of $ \mathbb{N}^0$
    we shall adhere to the notation introduced in section 2 of \cite{ER},
    except that the Young diagram of a composition or partition $ \mu $ shall be denoted
    $ \Y(\mu)$. We set 
$ {\mathcal Comp} := \bigcup_n \Comp $ and $ {\mathcal Par} := \bigcup_n \Par $. }}

\medskip
Suppose that $\mu\models n$. Then a $\mu$-\textit{tableau} is 
a bijection $\T: \mathbf{n} \to  \Y(\mu) $. 
If $\T $ is a $ \mu$-tableau we identify it with 
a labelling of the nodes of $ \Y(\mu) $, using the elements of $ \mathbf{n} $.
For example, if $\mu=(4,4,4,1)$ then
\begin{equation}\T = {\young(143\trece,26\diez8,57\once9,\doce)} \end{equation}
is the $\mu$-tableau $ \T$, that satisfies $ \T(1) = (1,1) , \T(4) = (1,2) , \T(3) = (1,3) $, etc.
If $\T $ is a $ \mu $-tableau we define $ Shape(\mathfrak{t}) :=   \mu $.
The set of $ \mu$-tableaux is denoted $ \Tab(\mu)$.
A $\mu$-tableau $\mathfrak{t}$ is \textit{row standard}
if the entries in $\mathfrak{t}$ increase from left to right in each
row and it is \textit{standard} if the entries also increase from
top to bottom {in each column}. The set of
row standard
$\lambda$-tableaux is denoted $\rstd(\lambda)$ and the set of standard $ \lambda $-
tableaux is denoted 
$\std(\lambda)$. For 
$\mu \models n $ we denote by $\mathfrak{t}^{\mu}$ (resp. $\mathfrak{t}_{\mu}   $) the standard
tableau in which the integers $1,2,\ldots,n$ are entered in
increasing order from left to right along the rows (resp. columns) of $\Y(\mu)$. For
example, if $\mu=(4,3)$
then $\T^{\mu}={\footnotesize\young(1234,567)}$ (resp.
$\T_{\mu}={\footnotesize\young(1357,246)}$).
If $ \lambda \in \Par $ and $ \T $ is a $ \lambda$-tableau we denote by $ \T^{\prime} $ the
\textit{conjugate} tableau. It is the $ \lambda^{\prime} $-tableau defined via
$ \T^{\prime}(i,j) := \T(j,i)$.

$\Si_n$ acts naturally on the right on $ \Tab(\lambda) $, via
composition of a bijection $ \T \in \Tab(\lambda) $ with an element $ \sigma \in \Si_n $
viewed as a bijection of $ {\bf n}$. 
For $\T \in \Tab(\lambda) $, we denote by
$d(\T)$ the unique element of $\Si_n$ such that
$\T=\T^{\lambda}d(\T)$.
We set $ w_{\lambda}:= d(\T_{\lambda}) \in \Si_n$. Then for any
$ \T \in \std(\lambda) $ 
we have that 
\begin{equation}\label{Murphyduality}
 d(\T) d(\T^{\prime})^{-1} =
w_{\lambda} \mbox{  and }  \ell(d(\T)) + \ell(d(\T^{\prime})) = \ell( w_{\lambda})
\end{equation}  
see for example the proof of Lemma 2.2 of \cite{M92}. 
The \textit{Young subgroup} {$\Si_{\lambda}$} associated with  $\lambda$
is the row stabilizer of $\mathfrak{t}^{\lambda}$.
It is a parabolic subgroup of $ \Si_n $ in the sense of Coxeter groups.
The set 
$  \left\{d(\T) \, | \,  \T  \in \rstd(\lambda) \right\}
$ 
is a set of \textit{distinguished right coset representatives} for $\Si_{\lambda}$ in $ \Si_n$,  that is
$\ell(w_0 d(\T))=\ell(w_0)+\ell(d(\T))$
for $w_0\in\Si_{\lambda}$. 
For any $ w \in \Si_n $ there is a unique decomposition
\begin{equation}\label{decompositionintialkindA}
  w = w_0 d(\T)\; \mbox{ where }\; w_0
  \in \Si_{  \lambda }\;\mbox{ and }\;  \T  \in  \rstd(\lambda)
\end{equation}
such that $\ell(w)=\ell(w_0)+\ell(d(\T))$, 
see Proposition 3.3 of \cite{Mathas} and chapter 1.10 of \cite{Hu}. The determination of
$ w_0 $ and $ \T $ can be realized in 
an algorithmic way that we now explain. 
Define first $ \s \in \Tab(\lambda) $ via $ w = d(\s)$. Let next
$ s_{j_1}, s_{j_2}, \ldots, s_{j_k}$ be a sequence of elements 
in $ S $ such that when setting $ \s_0 := \s, \, \s_1 := \s_0s_{j_1}, \, \s_2 := \s_1 s_{j_2},  \, \ldots, \,
\s_k := \s_{k-1}  s_{j_k} $ we have that $ j_i $ appears strictly above $ j_i+1$ in $ \s_j $ and that
$ \s_k $ is equal to $ \T^{\lambda}$ modulo a permutation of the row elements.
Such a sequence always exists.
Then 
$ w_0 = d(\s_k)$ and $ \T = \T^{\lambda} s_{j_k} s_{j_{k-1}} \dots s_{j_1} $, that is
$ d(\T) = s_{j_k} s_{j_{k-1}} \dots s_{j_1}$. Note that although $ w_0 $ and $ d(\T) $ are unique, the
sequence $ s_{j_1}, s_{j_2}, \ldots, s_{j_k}$ is in general not unique.
Here is an example of this algorithm, using 
$ \lambda = (2,3,2) $ and $ w = d(\s) $ given by 
\begin{equation}
\s= {\young(26,351,74)} \, .
\end{equation}
Then a possible sequence
of $ s_{j_i}$'s is given as follows
\begin{equation}\label{sequenceexample}
 \young(26,351,74) \, \stackrel{s_1} \longrightarrow \,  \young(16,352,74)
 \, \stackrel{s_4} \longrightarrow \,  \young(16,342,75)
 \, \stackrel{s_5} \longrightarrow \,  \young(15,342,76)
 \, \stackrel{s_4} \longrightarrow \,  \young(14,352,76)
 \, \stackrel{s_3} \longrightarrow \,  \young(13,452,76)   
\, \stackrel{s_2} \longrightarrow \,  \young(12,453,76)   
\end{equation}
and so $ d(\T) = s_2 s_3 s_4 s_5 s_4 s_1$
whereas $ w_0  = s_4 s_3 s_6 $ as can be seen from the last tableau of
(\ref{sequenceexample}). 
Note that 
\begin{equation}
\T = \T^{\lambda} s_2 s_3 s_4 s_5 s_4 s_1 =  \young(26,135,47) 
\end{equation}
which is also the tableau obtained from $ \s $ by ordering the rows. This holds in general.

\medskip
{\color{black}
  We write $\mu\lhd \nu$ for the usual dominance order on ${\mathcal Comp}$,
  see for example section 2 of \cite{ER} for
  the precise definition. For $ \lambda \in {\mathcal Comp}$ there is a natural extension of $ \lhd $ to
  a dominance order on 
  $ \rstd(\lambda) $, that we shall write the same way $  \lhd $,
  see once again \cite{ER} for the definition.}
We have that 
$\T^{\lambda} $ (resp. $ \T_{\lambda} $) is the unique maximal (resp. minimal) row standard $\lambda$-tableau
under the dominance order.

An \textit{$r$-multicomposition}, or just a multicomposition if confusion is not possible, 
of $n$ is an $r$-tuple
$\blambda=(\lambda^{(1)},\lambda^{(2)},\ldots,\lambda^{(r)})$ of
(possibly empty) compositions $\lambda^{(k)}$ such that
$\sum_{i=1}^{r}|\lambda^{(i)}|=n$. We call $\lambda^{(k)}$ the
$k$'th component of $\blambda$. An \textit{$r$-multipartition} is 
an $r$-multicomposition whose
components are all partitions.
The nodes $ \Y(\blambda) $ of an $r$-multicomposition $ \blambda $ are labelled by triples
$(x,y,p)$ with $p$ giving the component number and $ (x,y) $
the node of that component.
This is the
Young diagram for $ \blambda $ and is represented
graphically as the $ r$-tuple of Young diagrams of the components. For example, the
Young diagram of $\blambda=((3,3),(3,1),(1,1,1))$ is
$$\left( \yng(3,3) \;,\; \yng(3,1) \;,\; \yng(1,1,1)\; \right).$$
The set of $r$-multicompositions of $n$ is denoted by $ \MC  $ 
and the subset of $r$-multi\-partitions of $n$ is denoted by $ \MP  $.
If $\blambda \in \MC$, then 
a $\blambda$-\textit{multitableau} $ \bT $ is a
bijection $\bT: \mathbf{n} \to \Y({\blambda})$; it is 
identified with a labelling of $ \Y(\blambda) $ using the elements of
$\mathbf{n}$. The restriction of $\bT$ to $ \lambda^{(i)} $ is 
the $i$'th component $ \T^{(i)} $ of $\bT$ and we write $
\bT=(\T^{(1)}, \T^{(2)}, \ldots, \T^{(r)}) $. We say that $\bT$ is \textit{row
standard} if all its components are row standard, and
\textit{standard} if they are all standard. If $
\bT $ is a $ \blambda$-multitableau we write $ Shape(\bT) =
\blambda$. The set of all
$\blambda$-tableaux is
denoted by $\Tab(\blambda)$, the set of all
row standard $\blambda$-tableaux is
denoted by $\rstd(\blambda)$
and the set of all 
standard $\blambda$-multitableaux 
by $\std(\blambda)$.
In the following examples
\begin{equation}\label{11}
\bT=\left(\,\young(123,45)\;,\;\young(6,7,8)\;\right), \qquad
\Bs=\left(\,\young(278,14)\;,\;\young(56)\;,\;\young(3,9)
\;\right)
\end{equation}
$\bT$ is a standard multitableau whereas $\Bs$ is only a row standard multitableau. 
We denote by $\bT^{\blambda}$ (resp. $\bT_{\blambda}$) the
$\blambda$-multitableau in which $1,2,\ldots, n$ appear in order along the rows (resp. columns) of
the first component, then along the rows of the second component,
and so on. For example, in (\ref{11}) we have $ \bT= \bT^{\blambda}$ with 
$\blambda=((3,2),(1,1,1))$.
For each multicomposition $\blambda$ we define the Young subgroup
$\Si_{\blambda}$ as the row stabilizer of
$\bT^{\blambda}$. It is a parabolic subgroup in the sense of Coxeter groups.
For $\Bs$ a row standard $\blambda$-multitableau, we denote by
$d(\Bs)$ the unique element of $\Si_n$ such that
$\Bs=\bT^{\blambda}d(\Bs)$. 
The set 
$  \left\{d(\Bs) \, | \,  \Bs  \in \rstd(\blambda) \right\}
$ 
is a set of \textit{distinguished right coset representatives} for $\Si_{\blambda}$ in $ \Si_n$,  that is
$\ell(w_0d(\Bs))=\ell(w_0)+\ell(d(\Bs))$
for $w\in\Si_{\blambda}$.
Suppose that 
$ \blambda:= (\lambda^{(1)}, \ldots, \lambda^{(r)}) \in \MP$
and that $ \bT:= (\T^{(1)}, \ldots, \T^{(r)}) \in \Tab(\blambda)$. Then we define
the conjugation of $ \blambda $ and $ \bT $ componentwise, that is 
\begin{equation}
 \blambda^\prime:= ((\lambda^{(1)})^{\prime}, \ldots, (\lambda^{(r)})^{\prime}) \mbox{ and } 
 \bT^\prime:= ((\T^{(1)})^{\prime}, \ldots, (\T^{(r)})^{\prime}).
\end{equation}
    {\color{black}Note that this definition of conjugation differs from the definition that
      is often used in the literature, for example in \cite{HuMat}}.

Similarly to (\ref{decompositionintialkindA}), there is for any $ w \in \Si_n $ a unique decomposition
\begin{equation}\label{decompositionintialkindAA}
  w = w_0 d(\bT)\; \mbox{ where }\; w_0
  \in \Si_{  \blambda }\;\mbox{ and }\;  \bT  \in  \rstd(\blambda).
\end{equation}
Suppose that $ w= d(\Bs) $. Then 
applying the algorithm in (\ref{decompositionintialkindA})
on the tableau $ \s$ obtained from $ \Bs= (\s^{(1)}, \s^{(2)}, \ldots, \s^{(r)})$ by concatenating all its
components, with $\s^{(1)} $ on top followed by $  \s^{(2)} $ just below $\s^{(1)} $ and so on, 
we obtain $ w_0 $ and $ \bT $ in (\ref{decompositionintialkindAA}).

Suppose that $\blambda \in \MC$ and that $\bT \in \Tab(\blambda)$ is row standard.
Then for $j=1,\ldots,n$ we set 
$p_{\bT}(j):=k$ if $j$ appears in 
$\T^{(k)}$. We call $p_{\bT}(j) $ the
\textit{position} of $j$ in $\bT$.
When $\bT=\bT^{\blambda}$, we write $p_{\blambda}(j)$ for
$p_{\bT^{\blambda}}(j)$. We 
say that $\bT$ is of the
\textit{initial kind} if
$p_{\bT}(j)=p_{\blambda}(j)$ for all $j=1,\ldots,n$.

Suppose that $\blambda, \bmu \in \MC $.
We write $\blambda \unlhd \bmu$ if
$\lambda^{(i)}\unlhd \mu^{(i)}$ for all $i=1,\ldots,n$, this is the 
dominance order on $ \MC$. If $\Bs $ and $\bT$ are row standard
{multitableaux} and $m=1,\ldots,n$ we define $\Bs_{\mid \le m} $ and $ \bT_{\mid \le m}$
as for usual tableaux and
write $\Bs \unlhd \bT$ if $Shape(\Bs_{\mid \le m}) \unlhd
Shape(\bT_{\mid \le m}) $ for all $ m$.
Note that our dominance order $ \unlhd $ is different
from the dominance order on multicompositions and multitableaux that is
{\color{black}sometimes used} in the literature, for example in \cite{DJM}. 

To an $r$-multicomposition
$\blambda=(\lambda^{(1)},\ldots,\lambda^{(r)})$ we associate a
composition $ \|\blambda\| $ of length {$r$} in the following way
\begin{equation}\label{thecompositionas}
 \norm{\blambda}  := \left( | \lambda^{(1)}| ,\ldots,| \lambda^{(r)}|  \right).
\end{equation}
\color{black}\nomenclature[13]{$\norm{\blambda}$}{The composition associated with the multicomposition $\blambda$}Let $ \Si_{  \norm{\blambda} } $ be the associated Young subgroup.
Then $ w \in \Si_{  \norm{\blambda} } $ iff
$ \bT^{\blambda} w $ is of the initial kind.
For $ \bT = (\T^{(1)}, \ldots, \T^{(r)})  \in \Tab(\blambda)$
we define $  \|\bT\| \in \Tab( \| \blambda \|) $ via  
\begin{equation}
\| \bT\|  = (|\T^{(1)}|, \ldots, |\T^{(r)}|)
\end{equation}
where $ |\T^{(i)}| $ is the \textit{row reading} of $ \T^{(i)} $. 
For the multitableaux $ \bT $ and $ \Bs $ in (\ref{11}) we have for example
\begin{equation}
\|\bT\| = \young(12345,678)\, , \, \, \, \, \, \, \, \, \, 
\|\Bs\|= \young(27814,56,39).
\end{equation}
$ \Si_{  \norm{\blambda} } $ is yet another parabolic subgroup of $ \Si_n $, with 
corresponding distinguished right coset representatives
$ \left\{ \bT \in \Tab(\blambda) \;\middle|\;  \| \bT \| \in \rstd(\| \blambda \|) \right\} $.
Hence, for 
any $ w \in \Si_n $ there is a unique decomposition 
\begin{equation}\label{decompositionintialkind} w = w_0 d(\bT)\; \mbox{ where }\; w_0\in \Si_{  \norm{\blambda} }\;\mbox{ and }\; \bT \in 
  \Tab(\blambda)\;
\mbox{ with }\;
  \| \bT \| \in \rstd(\| \blambda \|).
\end{equation}
We define $ \Bs_0 \in \Tab(\blambda)$ via $ w_0 = \bT^{\blambda} d( \Bs_0) $; {it} is  of the initial kind.
We obtain $ \Bs_0 $ and $ \bT $ concretely by applying
the algorithm given in (\ref{sequenceexample}) with respect to $ w $ and $ \| \blambda \| $.
For example, if $ w = d(\Bs ) $ where
\begin{equation}
\Bs=\left(\,\young(2,6)\;,\;\young(3,51)\;,\;\young(7,4)
\;\right)
\end{equation}
then the calculations in (\ref{sequenceexample}) give us 
\begin{equation}
\Bs_0=\left(\,\young(1,2)\;,\;\young(4,53)\;,\;\young(7,6)
\;\right), \, \, \, \, \, \, \, \, \, 
\bT=\left(\,\young(2,6)\;,\;\young(1,35)\;,\;\young(4,7)
\;\right)
\end{equation}
where we used $ w_0 = s_4 s_3 s_6 $ to determine $ \Bs_0$. In general, 
it follows from the algorithm that $ \Bs_0 $ is row standard if $ \Bs $ is row standard.

Suppose that $ \Bs \in \rstd(\blambda) $ and $ 
\Bs_1 \in \rstd(\blambda_1) $ 
with $ \norm{\blambda} = \norm{\blambda_1} $.
Let $ w = d(\Bs),  w_1 = d(\Bs_1 ) $
and let $ w= w_0 d(\bT), w_1= (w_1)_0 d(\bT_1) $ be the decompositions of $ w $ and $ w_1 $ 
as in (\ref{decompositionintialkind}). 
Suppose furthermore that 
$d(\bT) = d(\bT_1) $. 
Then we have that 
\begin{equation}\label{thenweave}
\Bs \unlhd \Bs_1 \mbox{ if and only if } \Bs_0 \unlhd (\Bs_1)_0.
\end{equation}

For $ \blambda $ any $r$-multi\-partition 
we define $ w_{\blambda} \in \Si_n $ as $ w_{\blambda} := d(\bT_{\blambda}) $.
Suppose now that $ \Bs \in \Tab(\blambda)  $ and let $ \Bs^{\prime} $ be the
conjugate multitableau.
Set
$ w = d(\Bs) $ and $ w^{\prime} = d(\Bs^{\prime}) $ with 
decompositions 
$ w = w_0 d(\bT)  $ and $ w^{\prime} = (w^{\prime})_0 d(\bT^{\prime})  $ 
as in (\ref{decompositionintialkind}). Then it follows from the algorithm described above that
$ d(\bT) = d(\bT^{\prime}) $ and that $ (w^{\prime})_0 = d(\bT_0^{\prime}) $ where
$ w_0  = d(\bT_0) $. Moreover, since $ \bT_0 $ is of the initial kind 
we get via (\ref{Murphyduality}) that
\begin{equation}\label{22}
 d(\Bs) d(\Bs^{\prime})^{-1} =
  w_{\blambda} \mbox{  and }  \ell(d(\bT_0)) + \ell(d(\bT_0^{\prime})) = \ell( w_{\blambda}). 
\end{equation}  
Note that we cannot here replace $ d(\bT_0) $ and $ d(\bT_0^{\prime}) $ by
$ d(\Bs) $ and $ d(\Bs^{\prime}) $.

\medskip

Let us now recall the $bt$-algebra of the title of the paper.
It was originally introduced by Aicardi and Juyumaya in \cite{AJ1}.
\begin{defi}\label{braidsties}
  Let $n$ be a positive integer.  The braids and ties algebra
  ${\cal E}_n = \E $, or the $bt$-algebra, is the associative $S$-algebra generated by the
elements $g_1,\ldots,g_{n-1} $ and $e_1, \ldots, e_{n-1}$, subject to the
following relations:
\begin{alignat}{3}
g_ig_j&=g_jg_i&&\quad\mbox{ for } |i-j|>1\label{E1}\\
g_ie_i&=e_ig_i&& \quad \mbox{ for all } i \label{E2}\\
g_ig_{j}g_i&=g_{j}g_ig_{j}&&\quad\mbox{ for }  |i-j|=1 \label{E3}\\
e_ig_{j}g_i&=g_{j}g_ie_{j}&&\quad\mbox{ for }  |i-j|=1 \label{E4}\\
e_ie_{j}g_j&=e_{i}g_je_{i}=g_{j}e_ie_{j}   &&\quad\mbox{ for }  |i-j|=1 \label{E5}\\
e_ie_j&=e_je_i&&\quad\mbox{ for all } i,j\label{E6}\\
g_ie_j&=e_jg_i&&\quad\mbox{ for } |i-j|>1\label{E7}\\
e_i^2&=e_i&&\quad\mbox{ for all } i\label{E8} \\
g_i^2&=1+(q-q^{-1})e_ig_i  && \, \, \, \, \, \, \, \, \,  \mbox{  for all }\; i\label{E9}.
\end{alignat}
\end{defi}

There is a diagrammatic interpretation of these relations that explains the name of $ \E$
and that we now briefly indicate even though it only plays a minor role in the present paper.
It uses $n$ strands connecting a northern
and a southern border. Under the interpretation, the $ g_i$'s correspond 
to usual simple \textit{braids} involving the $i$'th and the $i+1$'st strand{\color{black}s} 
and the $ e_i$'s correspond to \textit{ties} involving the $ i$'th and $i+1$'st strands. 
The relations (\ref{E4}) and (\ref{E5}) are for example visualized as follows.
\begin{equation}\label{relacion1}
  \raisebox{-.5\height}{\includegraphics[scale=0.4]{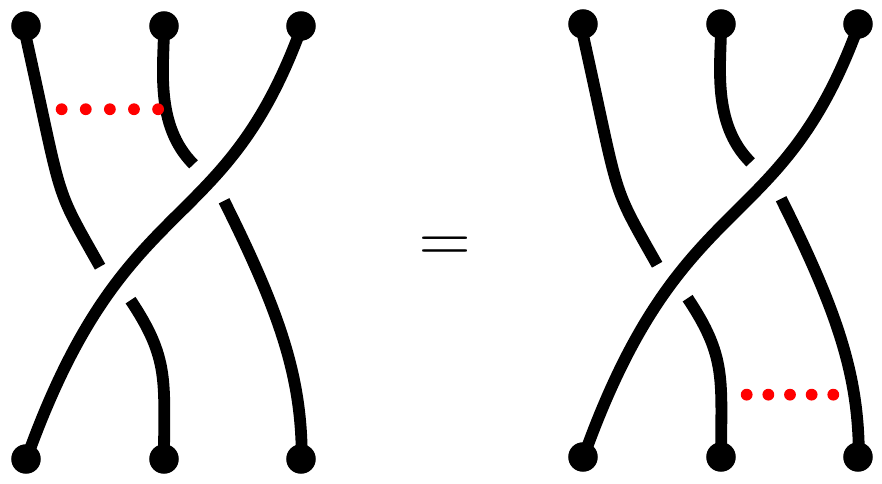}} \, \, \,\, \, \, \, \,\, \, 
\end{equation}
\begin{equation}\label{relacion2}
  \raisebox{-.5\height}{\includegraphics[scale=0.4]{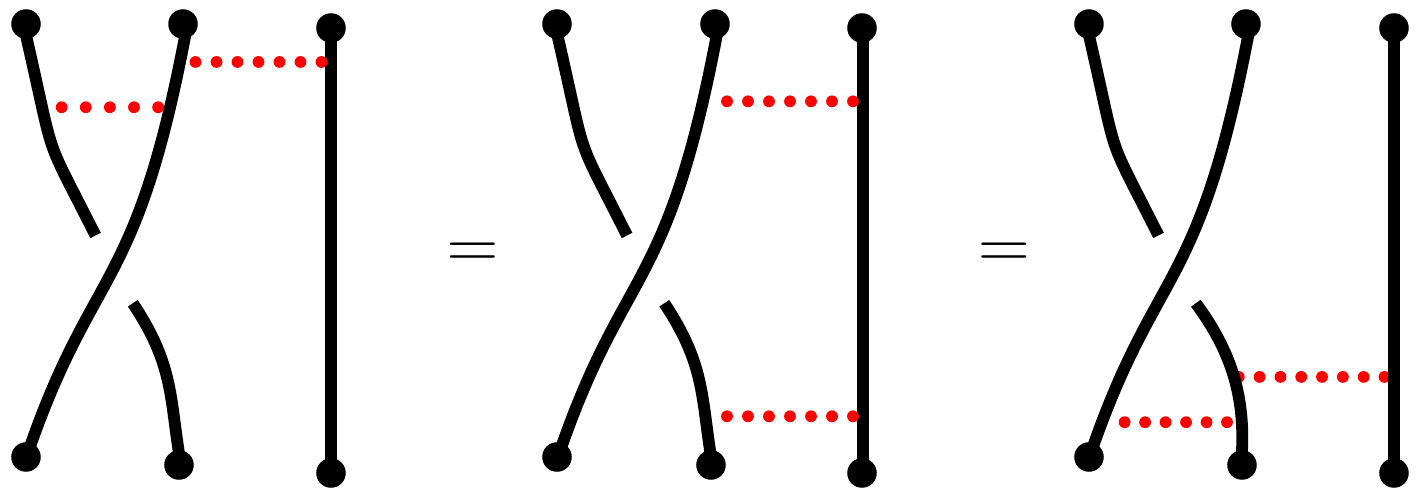}} \, \, \,\, \, \, \, \,\, \, 
\end{equation}

\medskip
If $ \kk $ is a commutative ring containing an invertible element, that we also denote $q  $,
then we define the specialized algebra ${ \mathcal E}^{\kk}_n(q) $
via ${ \mathcal E}^{\kk}_n(q)  := \E \otimes_S \kk $ 
where $ \kk $ is made into an $S$-algebra by mapping $
q \in S $ to $ q \in \kk $.

\medskip
The concept of \textit{set partitions} is important in the study of $ \E$.
Recall that a set partition $A=\{I_1,I_2,\ldots I_r\}$ of $\mathbf{n}$ is a set of 
nonempty and disjoint subsets 
of $\mathbf{n}$ 
whose union is $\mathbf{n}$, in order words, $ A $ is the set of classes of
an equivalence relation on $\mathbf{n}$.
We refer to the $ I_i$'s as
the \textit{blocks} of $A$. We denote by $\mathcal{SP}_n$ the set of all set partitions of
$\mathbf{n}$. 
There is a canonical poset structure on $
{\cal SP}_n$ defined as follows. Suppose that $A=\{I_1,I_2,\ldots,
I_k\},  B=\{J_1,I_2,\ldots, J_l \} \in {\cal
SP}_n$. Then we say that $ A \subseteq B $ if each $ J_j $ is a
union of some of the $ I_i$'s.

\medskip
Returning to $ \E $, for 
$ 1 \le i <j  \le  n $ we define $ E_{ij} \in \E$ by $ E_{ij} = e_i $ if $  i= j-1$, and
recursively downwards on $ i $ via 
\begin{equation}\label{makessense}
  E_{ij}:=g_{i}  E_{i+1,j} g_{i}^{-1}.
\end{equation}
Note that $ g_{i} $ is invertible with inverse $ g_i^{-1} = g_i +(q^{-1} -q) e_i $, as follows from 
(\ref{E9}), so that (\ref{makessense}) makes sense.
For $ A \in {\cal SP}_n $ we define $ E_A \in \E$ as 
\begin{equation}\label{runsoverpairs}
E_A := \prod_{i,j}   E_{ij}
\end{equation}
where the product runs over pairs $ (i,j),  $ such that $ i < j $ and such that $i $ and $ j $
belong to the same block of $ A$. One checks that the $E_{ij}$'s commute and so the
product is independent of the order in which it is taken.
For $ w = s_{i_1} \cdots s_{i_k} \in \Si_n$ a reduced expression for $w$ we define
$ g_w := g_{i_1} \cdots g_{i_k} \in \E $. By the relations it is independent of the choice of
reduced expression.
We have the following relations in $ \E$
\begin{equation}{\label{wehavetherelations}}
E_A g_w = g_w E_{ A w}  \mbox{ and } E_A E_B = E_C \mbox{  for } w\in \Si_n , A,B  \in \mathcal{SP}_n
\end{equation}
where $ C \in \mathcal{SP}_n $ is minimal with respect to $ A \subseteq C, B \subseteq C$.
Moreover, it was shown in \cite{Ry} that $ \E $ has an $S$-basis of the
form $ \left\{ E_{A} g_w | A \in \mathcal{SP}_n, w \in \Si_n \right\} $ and so, in particular,
the dimension of $ \E $ is $ n! b_n $ where $ b_n $ is the Bell number,
that is the cardinality of $ \mathcal{SP}_n $.



\section{Cellular basis for $\E$.}
In the paper \cite{ER}, a cellular basis $ \{ m_{\es\et} \} $ for $ \E $ was constructed.
Moreover, it was shown that $ \{ m_{\es\et} \} $ induces cellular bases for
interesting subalgebras of $ \E$ and via this an isomorphism between $ \E $ and a direct sum of matrix algebras
over certain wreath products of Hecke algebras was established.

\medskip
In this section we first recall the construction of $ \{ m_{\es\et} \} $ and at the same time
we introduce a new, \textit{dual}, cellular basis for $ \E$, that we denote $ \{ n_{\es\et} \} $.
We next prove,  in a series of Lemmas, certain compatibility properties between the 
dominance order associated with $ \E $ and the multiplication of elements from the
two bases. These compatibility properties are known in the Hecke algebra setting, and
our proofs are essentially reductions, although not completely trivial, to that setting.

\medskip
Let us first recall the formal definition of cellular algebras, as introduced in \cite{GL}.
\begin{defi}{\label{cellular}} Let $\mathcal{R}$ be an integral domain. Suppose that $A$ is an
  $\mathcal{R}$-algebra. Suppose that $(\PP,\leq)$ is a poset
  and that for each $\lambda\in \PP$ there is a
finite indexing set $T(\lambda)$ (the '$\lambda$-tableaux') and elements $c_{\s\T}^{\lambda}\in
A$ such that
\begin{equation}
  \mathcal{C}=\left\{c_{\s\T}^{\lambda}\mid \lambda\in \PP \mbox{ and } \s,\T\in T(\lambda)\right\}
\end{equation}  
is an $\mathcal{R}$-basis of $A$. The pair $(\mathcal{C},\PP)$ is a
\textit{cellular basis} for $A$ if
\begin{enumerate}\renewcommand{\labelenumi}{\textbf{(\roman{enumi})}}
\item The $\mathcal{R}$-linear map $*:A\to A$ determined by
$(\mathop{c_{\s\T}^{\lambda}})^*=c_{\T\s}^{\lambda}$ for all
$\lambda\in\PP$ and all $\s,\T\in T(\lambda)$ is an algebra
anti-automorphism of $A$.

\item For any $\lambda\in \PP,\; \T\in T(\lambda)$ and $a\in A$
there exist elements $r_{\V}\in \mathcal{R}$ such that for all $\s\in T(\lambda)$
$$c_{\s\T}^{\lambda} a \equiv \sum_{\V\in T(\lambda)} r_{\V}c_{\s\V}^{\lambda} \mod{A^{\lambda}}$$
where $A^\lambda$ is the $\mathcal{R}$-submodule of $A$ with basis
$\left\{c_{\U\V}^{\mu}\mid \mu\in \PP,\mu>\lambda \mbox{ and }
\U,\V\in T(\mu)\right\}$.
\end{enumerate}
If $A$ has a cellular basis we say that $A$ is a \textit{cellular
algebra} and $(\PP,T,\mathcal{C},\ast)$ is called the \textit{cell datum} for $A$.
\end{defi}

\medskip
We now recall the, somewhat lengthy, construction of the cell datum for $ \E$. 
For more details on this construction, {\color{black} and for all unexplained notions}, the reader should consult \cite{ER}.

We first consider a certain ideal decomposition of $ \E$.
Let $  \mu_{ {\cal SP}_n}  $ be the M\"o
bius function for the
  lattice $ ({\cal SP}_n, \subseteq)$, {\color{black} see for example section 6.1 of \cite{ER} for 
    precise formulas and examples}.
Using it, we define 
for $A\in \mathcal{SP}_n$ orthogonal idempotents $\mathbb{E}_A \in \E$
via 
\begin{equation}{\label{moebius}}
\mathbb{E}_A:= \sum_{B: A \subseteq B }\mu(A,B)E_B.
\end{equation}
The following is Proposition 39 in \cite{ER}.
\begin{propos}\label{commuting} The following properties hold.
\begin{itemize}
\setlength\itemsep{-1.3em}
\item[(1)] $ \{ \mathbb{E}_A | A \in \mathcal{SP}_n \} $ is a set of orthogonal idempotents of $ \E$. \newline
\item[(2)] For all $w\in\Si_n$ and $A\in \mathcal{SP}_n$ we have $\be_A g_w=g_w\be_{Aw}$. \newline
\item[(3)] For all $A\in \mathcal{SP}_n$ we have $\be_{A}E_B=\left\{\begin{array}{ll}\be_{A}&\mbox{ if } B\subseteq A\\
0&\mbox{ if } B\not\subseteq A.  \end{array}\right.$
\end{itemize}
\end{propos}

Let $ \alpha \in  \Par$.
A set partition
$A=\{I_{1},\ldots,I_{k}\}$ of $\mathbf{n}$ is said to be of
\textit{type} $\alpha $ if there exists a permutation $\sigma$ such
that $(|I_{i_{1\sigma}}|,\ldots,|I_{i_{k\sigma}}|)=\alpha$. 
We write $|A|=\alpha$ if $A\in \mathcal{SP}_n$ is of type $\alpha$
and we let $ \mathcal{SP}_n^{\alpha} $ be the set of set partitions
of type $ \alpha$. 
For $\alpha\in \Par$  we define the following element $\be_{\alpha}$ of $\E$
\begin{equation}{\label{Ealpha}}
\be_{\alpha}:=\sum_{{A\in \mathcal{SP}_n},\,  |A|=\alpha}\be_{A}.
\end{equation}
Then the $ \be_\alpha $'s form a family of central
orthogonal idempotents in $\E$ such that $ \sum_{ \alpha \in \Par} \be_\alpha = 1 $.
As a consequence we have the following decomposition of $\E$ into a direct sum of two-sided ideals
\begin{equation}\label{descE}
\E=\bigoplus_{\alpha\in\Par} \Ea
\end{equation}
where $\Ea:=\be_{\alpha}\E$.
{Each ideal $ \Ea $ is an $S$-algebra with identity $
\be_{\alpha} $}
and the set 
\begin{equation}\label{basessubalgebras}
 \{\be_{A}g_w\mid w\in\Si_n,\,|A|=\alpha\}
\end{equation}
is an $S$-basis for $\Ea$.
For $ \kk $ an arbitrary field {\color{black}containing} a nonzero element $ q $ we have a specialized version of $ \Ea $
\begin{equation}
\Eak := \Ea \otimes_{S} \kk.
\end{equation}

In view of the decomposition in (\ref{descE}), in order to give the cell datum for $ \E $ it is enough to 
give the cell datum for each $\Ea$.


\medskip
We now explain the poset denoted $ \PP $ in Definition
\ref{cellular}. We first describe $ \PP $ as a set.
For $ \blambda = (\lambda^{(1)}, \ldots , \lambda^{(r)}) \in \MC $
we define sets
$ I_i  $ via
 $ I_1:=\{1, 2, \ldots, |\lambda^{(1)}|\} $, 
$ I_2:=\{|\lambda^{(1)}|+1, |\lambda^{(1)}|+2, \ldots, |\lambda^{(1)}|+|\lambda^{(2)}|\} $, and so on.
Leaving out the empty $ I_i $'s we obtain a set partition
in $ \mathcal{SP}_n $ that we denote
$ A_{\blambda} $.
For $ \alpha \in \Par $ we say that $ \blambda $ is of type $ \alpha $ if $ A_{\blambda} $ is of
type $ \alpha$.

An $ r$-multipartition $ \blambda
= (\lambda^{(1)}, \ldots , \lambda^{(r)}) $ is said to be \textit{increasing}
if $\lambda^{(i)} < \lambda^{(j)} $ iff $ i < j $, where  $ < $ is
any fixed extension of the dominance order $ \lhd $ to a
total order.

\medskip
For $ \alpha \in \Par$, we now define 
$ {\mathcal L}_n(\alpha) $ to be 
the set of pairs $ \Lambda =
(\blambda \mid \bmu) $ where $ \blambda = ( \lambda^{(1)}, \ldots ,
\lambda^{(r)})$ is an increasing $r$-multipartition of $ n $ of type $ \alpha ${\color{black}.}
Let $ \lambda^{(i_1)}  < \cdots < 
\lambda^{(i_s)} $ be the distinct $ \lambda^{(i)} $'s and let $ m_j $ be the multiplicity of
$ \lambda^{(i_j)} $ in $ ( \lambda^{(1)}, \ldots ,\lambda^{(r)}) $. Then we require that $ \bmu $ be
an $s$-multipartition $ \bmu = (\mu^{(1)}, \ldots,  \mu^{(s)}) $ of $ r $
where $ \mu^{(j)} \vdash m_j $. For $ \PP$ we choose $ {\mathcal L}_n(\alpha) $.

\medskip

We next describe the poset structure on $ {\mathcal L}_n(\alpha) $. 
Suppose that $ \Lambda =(\blambda \mid \bmu) $ and $
\overline{\Lambda} =(\overline{\blambda}\mid \overline{\bmu} ) $ are
elements of $  {\mathcal L}_n(\alpha) $, such
that $ \blambda = (\lambda^{(1)}, \ldots, \lambda^{(r)}) $ and
$ \overline{\blambda} = (\overline{\lambda^{(1)}}, \ldots, \overline{\lambda^{(r)}}) $.
Then we write $ \blambda \lhd_1 \overline{\blambda} $ if
there exists a permutation $ \sigma $ such that
$ (\lambda^{ (1 \sigma) }, \ldots, \lambda^{ (r\sigma )}) \lhd (\overline{\lambda^{(1)}}, \ldots,
\overline{\lambda^{(r)}})$
where $\lhd$ is the dominance order on $r$-multi\-partitions, introduced above.
We then write of $ \Lambda \lhd \overline{\Lambda} $ if $ \blambda \lhd_1 \overline{\blambda} $
or if $ \blambda =\overline{\blambda} $ and
{$ \bmu \lhd \overline{\bmu}$.}
As usual we set $ \Lambda \unlhd \overline{\Lambda} $ if $ \Lambda \lhd \overline{\Lambda} $ or if
$ \Lambda = \overline{\Lambda} $.
This is the description of $ {\mathcal L}_n(\alpha) $ as a poset.

\medskip


For $ \Lambda=(\blambda \mid  \bmu) \in { \mathcal L}_n(\alpha) $ as above,
we now define the concept of $ \Lambda $-tableaux.
Suppose that $ \et $ is a pair $\et=(\bT \mid \bu) $.
Then $ \et $ is called a $ \Lambda $-\textit{tableau} if $ \bT $ is a
$\blambda$-multitableau and 
$ \bu $ is a {$\bmu$-multitableau of the initial kind.}
If $ \et $ is a $ \Lambda $-tableau we define $ Shape(\et) := \Lambda$.
We let $ \Tab(\Lambda) $ denote the set of $ \Lambda$-tableaux.
We 
say that $ \et = (\bT \mid \bu) \in \Tab(\Lambda) $ is row standard
if its ingredients are row standard multitableaux in the 
sense of the previous section
and if moreover $\bT$ is an \textit{increasing}
multitableau.
By increasing we here mean
that whenever $ \lambda^{(i)} = \lambda^{(j)} $ we have that $ i < j $ if and only if
$ \min ( \T^{
  (i)}) < \min (\T^{(j)}) $ where $ \min (t) $ is the function that reads off the minimal entry of the tableau $ t $.
We say that $ \et = (\bT \mid \bu) \in \Tab(\Lambda) $ is standard if it is row standard 
and if its ingredients $ \bT $ and $ \bu$ are standard multitableaux and
we define $ {\rm Std}(\Lambda) $ to be the set of all standard $\Lambda$-tableaux.

Note that our notation is here deviating slightly from the one used in \cite{ER}
where the condition on $\bT $ to be increasing was required for standardness of $\et $, but not for
row standardness of $\et $.
In the following examples 
\begin{equation}\begin{array}{c}
\es :={\footnotesize\left(\;\left(\;\young(1,8)\;,\young(56)\;,\young(39)\;,\young(24,7)\;\right) \big|
  \left(  \young(1) \, , \, \young(2,3)\; ,\young(4) \right)\;\right)}
\\[4mm]
     \et :={\footnotesize\left(\;\left(\;\young(1,8)\;,\young(35)\;,\young(69)\;,\young(24,7)\;\right)  \bigm|
\left(  \young(1) \, , \young(2,3) \, , \,\young(4)\;\right)\;\right)}
\end{array}\end{equation}
$ \es $ is not row standard but $ \et $ is.


\medskip
We fix the following combinatorial notation.
Suppose that $ \Lambda=(\blambda \mid \bmu) \in  {\cal L}_n(\alpha) $ 
with
$ \blambda = ( \lambda^{ (1 ) }, \ldots,  \lambda^{ (r ) }) $ and
$ \bmu = ( \mu^{ (1 ) }, \ldots,  \mu^{ (s ) })$.  
By construction, the numbers $ m_i := |   \mu^{ (i ) }  | $ are the 
multiplicities of equal $ \lambda^{(i)}$'s.
Let $\Si_{\Lambda} \leq \Si_{n}$ be the stabilizer subgroup of the set
partition $ A_{\blambda} {=\{I_1, I_2, \ldots,
I_m\}} $. Then
the multiplicities give rise to the subgroup $ \Si^m_{\Lambda} $ of $ \Si_{\Lambda}$
consisting of the order preserving permutations of
the blocks of $ A_{\blambda} $ that correspond to equal $ \lambda^{(i)}$'s.
We have $ \Si^m_{\Lambda} \leq  \Si_{{\Lambda}}$, and 
as an abstract group 
\begin{equation}\label{abstractgroup}  \Si^m_{\Lambda}  = \Si_{m_1}  \times \cdots \times \Si_{m_s}.
\end{equation}  
For example, if $ \blambda= ( (1^2), (1^2), (1^2), (2), (2), (1^3), (1^3), (1^3), (2,1) ) $
then 
\begin{equation*}\label{aboveexample}
\begin{array}{c}
\bT^{\blambda}={\footnotesize\left(\;\young(1,2)\,,\,\young(3,4)\,,\,\young(5,6)\,,\,\young(78)\,,\,\young(9\diez)
  \,,\,\ytableausetup{centertableaux,boxsize=1.3em}\begin{ytableau}
    11\\12\\13\end{ytableau}\,,\,\begin{ytableau}14\\15\\16
    \end{ytableau}\,,\,\begin{ytableau}17\\18\\19
    \end{ytableau}\;\,,\,\begin{ytableau}20 &21\\22
    \end{ytableau}\; 
    \right)} \,\mbox{ and } \\
A_{\blambda} = \left((1,2) , (3,4) , (5,6) , (7,8) , (9,10) , (11,12,13) , (14,15,16) ,
(17,18,19) , (20,21,22) \right) 
\end{array}
\end{equation*}
and so $  \Si^m_{\Lambda}  = \Si_{3}  \times \Si_{2}\times \Si_{3} \times  \Si_{1}$.
For $ I_i$ and $ I_{i+1}$ two consecutive blocks of $ A_{\blambda} $ such
that $ \lambda^{(i)} = \lambda^{(i+1)} $ we let 
$ B_i \in \Si_{n}$ be the element of 
minimal length, in the sense of Coxeter groups, 
that interchanges $ I_i $ and $ I_{i+1} $. For example in the above example
we have $ B_1= (1,3)(2,4), B_2 =(3,5)(4,6) $, and so on.
$  \Si_{\Lambda} $ is generated by the $ B_i$'s such that $ \lambda^{(i)} = \lambda^{(i+1)} $, 
in fact it is 
a Coxeter groups on these 
$B_i$'s. 

\medskip
It is important that the group algebra $ S \Si^{m}_{\Lambda} $
can be realized as a subalgebra
of $ \E$.
Let $ \BB_i \in \Ea $ be the element introduced in equation (113) of \cite{ER}.  
{\color{black}Then} the following is (1) of Lemma {\color{black} 46} of \cite{ER}.

\begin{lem}{\label{embedding}}
Suppose that $ \Lambda=(\blambda \mid \bmu) \in  {\mathcal
L}_n(\alpha) $. Then we have an $S$-algebra embedding
$\iota:S\Si^m_{\Lambda} \hookrightarrow \Ea $ given by $ B_i \mapsto \BB_i $
where $ \lambda^{(i)} = \lambda^{(i+1)} $.
\end{lem}

\medskip
We are now ready to recall the construction of the basis elements $ \{ m_{\es\et} \} $ of the cell datum for
$ \Ea$, and at the same time we
introduce the new basis elements $ \{ n_{\es\et} \} $
for $ \Ea$.
For each $ \Lambda \in {\cal L}_n(\alpha) $, we first define 
elements $ m_{\Lambda} $ and $ n_{\Lambda} $ that act as starting points for the bases.
Suppose that $ \Lambda = (\blambda \mid \bmu ) $ is as above
with $ \blambda = (\lambda^{(1)}, \ldots, \lambda^{(r)}) $ and $ \bmu = (\mu^{(1)}, \ldots,\mu^{(s)})$.
We then define $ m_{\Lambda}  \in  \Ea $ and $ n_{\Lambda} \in \Ea$ as follows
\begin{equation}{\label{startingbasiselement}}
  m_{\Lambda}:= {\be}_{\blambda}   x_{\blambda}     b_{\bmu}, \, \, \, \, \, \, 
  n_{\Lambda}:= {\be}_{\blambda}   y_{\blambda}     c_{\bmu}.
\end{equation}
Let us explain the factors of these products.
The idempotent ${\be}_{\blambda} $ has already been introduced.
The factors $ x_{\blambda} \in \Ea $ and $ y_{\blambda} \in \Ea $
are defined as
\begin{equation}\label{secondlyfirst} x_{\blambda} :=  \sum_{w\in\Si_{\blambda}}q^{\ell(w)}g_{w}
\, \, \, \, \, \, \, \, 
\mbox{   and     }
\, \, \, \, \, \, \, \, 
y_{\blambda} :=  \sum_{w\in\Si_{\blambda}}(-q)^{-\ell(w)}g_{w}.
\end{equation}
Note that $ {\be}_{\blambda} $ commutes with $ x_{\blambda} $ and $ y_{\blambda} $ and that 
for $  w \in \Si_{\blambda} $, we have that
\begin{equation}{\label{justasinYH}}
{\be}_{\blambda} x_{\blambda} g_w = g_w  {\be}_{\blambda} x_{\blambda} = q^{l(w)} {\be}_{\blambda} x_{\blambda} \, \, \, \, \mbox{  and     } \, \, \, \,
{\be}_{\blambda} y_{\blambda} g_w = g_w  {\be}_{\blambda} y_{\blambda}
= (-q)^{-l(w)} {\be}_{\blambda} y_{\blambda}.
\end{equation}
In the Hecke algebra case, corresponding to $ r=1$, the elements defined in (\ref{secondlyfirst})
are the elements denoted 
$  x_{\lambda \lambda} $ and $  y_{\lambda \lambda} $ in \cite{Mur95}.

In order to define the factors $ b_{\bmu} $ and $ c_{\bmu} $ we recall from 
(\ref{abstractgroup}) the decomposition
$\Si^m_{\Lambda}  = \Si_{m_1}  \times \cdots \times \Si_{m_s} $. 
Let $ x_{\bmu}(1) \in \Si^m_{\Lambda}  $ and $ y_{\bmu}(1) \in \Si^m_{\Lambda} $ 
be the $ q= 1 $ specializations of the elements given in (\ref{secondly})
corresponding to the 
multipartition $\bmu$, that is
\begin{equation}\label{secondly} x_{\bmu}(1) :=  \sum_{w\in\Si_{\bmu}}g_{w} 
\, \, \, \, \, \, \, \, 
\mbox{   and     }
\, \, \, \, \, \, \, \, 
y_{\bmu}(1) :=  \sum_{w\in\Si_{\bmu}}(-1)^{-\ell(w)}g_{w}.
\end{equation}
Then $   b_{\bmu} \in \Ea$ and $   c_{\bmu} \in \Ea $ are defined via 
\begin{equation}
  b_{\bmu}:=\iota(x_{\bmu}(1)) \, \, \, \, \mbox{  and     } \, \, \, \,
  c_{\bmu}:=\iota(y_{\bmu}(1))
\end{equation}
where $ \iota:S\Si^m_{\Lambda} \hookrightarrow \Ea $ is the embedding from Lemma \ref{embedding}.

\medskip
Let $ \et^{\Lambda} $ (resp. $ \et_{\Lambda} $) be the $ \Lambda$-tableau given 
by $  \et^{\Lambda} :=(\bT^{\blambda} \mid \bT^{\bmu}) $
(resp. $  \et_{\Lambda} :=(\bT_{\blambda} \mid \bT_{\bmu}) $). 
Let now $ \es =(\Bs \mid  \bu) $ be an arbitrary 
$\Lambda$-tableau. 
Then 
we set $ d(\es) := (d({\Bs}) \mid \iota(d(\bu))) $ 
where $ d(\Bs) \in \Si_n $ satisfies 
$ \bT^{ \blambda } d(\Bs) =  \Bs $ and $ d(\bu) \in \Si^m_{\Lambda}$
satisfies $ \bT^{ \bmu} d(\bu) = \bu $. For simplicity, we write $(d({\Bs}) \! \!\mid d(\bu)) $
instead of $ (d({\Bs}) \mid \iota(d(\bu))) $.
Note that since
$\bu=(\U_1, \ldots,\U_s)$ is of the initial kind, we have a decomposition
$ d(\bu) = (d(\U_1),  \ldots, d(\U_s)) $, according to the decomposition in 
(\ref{abstractgroup}) and so 
\begin{equation}
  \BB_{d(\bu)}=  \BB_{d(\U_1)}\cdots \BB_{d(\U_s)}.
\end{equation}  
{\color{black}We can now finally recall the construction of the cellular basis $ m_{\es\et}$ from \cite{ER} and
at the same time introduce the new {\it dual cellular basis} $ n_{\es\et} $ for $ \Ea $.}
For $\es=(\Bs \mid \bu),\, \et=(\bT \mid \bv) $ row standard $\Lambda$-tableaux we define
\begin{equation}{\label{maindefM}}
m_{\es\et}:=   { \be}^{}_{\blambda} g_{d(\Bs)}^{\ast}  x^{}_{\blambda}
\BB_{d(\bu)}^{\ast}   b^{}_{\bmu}\BB^{}_{d(\bv)}  \,
g^{}_{d(\bT)}
\end{equation}
\begin{equation}{\label{maindefN}}
n_{\es\et}:=  { \be}^{}_{\blambda} g_{d(\Bs)}^{\ast}  y^{}_{\blambda}
\BB_{d(\bu)}^{\ast}  c^{}_{\bmu}\BB^{}_{d(\bv)}  \,
g^{}_{d(\bT)}.
\end{equation}

\medskip
\noindent
We set $  b^{}_{\bu \bv} :=  \BB_{d(\bu)}^{\ast}   b^{}_{\bmu}\BB^{}_{d(\bv)} $ 
and $  c^{}_{\bu \bv} :=  \BB_{d(\bu)}^{\ast}   c^{}_{\bmu}\BB^{}_{d(\bv)} $ and have then also
\begin{equation}
m_{\es\et}=   { \be}^{}_{\blambda} g_{d(\Bs)}^{\ast}  x^{}_{\blambda} b^{}_{\bu \bv} \,
g^{}_{d(\bT)}
\end{equation}
\begin{equation}
n_{\es\et}=   { \be}^{}_{\blambda} g_{d(\Bs)}^{\ast}  y^{}_{\blambda} c^{}_{\bu \bv} \,
g^{}_{d(\bT)}.
\end{equation}
In the case where $ \es = \et^{\Lambda} $ we write for simplicity
\begin{equation}\label{forsimplicity}
m_{\et} :=   m_{\et^{\Lambda} \et}  \mbox{  and  } 
n_{\et} :=   n_{\et^{\Lambda} \et}.
\end{equation}

In \cite{ER} the following Theorem was proved. 
\begin{teo}
$ \Ea $ is a cellular basis with cell datum $ ({\cal L}_n(\alpha), \std(\Lambda), m_{\es\et}, \ast)$.
\end{teo}  

Making the straightforward adaptations of the proofs given in \cite{ER}, we also have
the following Theorem.
\begin{teo}
$ \Ea $ is a cellular basis with cell datum $ ({\cal L}_n(\alpha), \std(\Lambda), n_{\es\et}, \ast)$.
\end{teo}  

By general cellular algebra theory, specialization induces a cell datum on 
$ \Eak$ as well. We shall denote it the same way; in particular $ m_{\es\et} $ and $  n_{\es\et} $
may refer to elements
of both $ \Ea $ and $ \Eak$.

The definition of cellular basis does not stipulate a partial order on the tableaux, but
in our setting there is a natural such {\color{black} \textit{dominance}} order $\color{black}{ \lhd}$
on tableaux
that we shall need. 
It was already introduced in \cite{ER}, {\color{black}{see the paragraph below {\bf Remark 43}
    of {\it loc. cit.} for
    the precise definition.}}
As usual we set $ \et \unlhd \overline{\et} $ if $ \et \lhd \overline{\et} $ or
$ \et = \overline{\et} $.

\medskip
We now introduce the concept of \textit{conjugation} of the elements
of ${\cal L}_n(\alpha) $ and of their tableaux; this concept is also not present in general cellular algebras.
{\color{black}{Note that it was not formally defined in \cite{ER}}}.
Suppose that $ \Lambda = (\blambda \mid \bmu ) \in {\cal L}_n(\alpha)$,  
with $ \blambda = (\lambda^{(1)}, \ldots, \lambda^{(r)}) $ and $ \bmu = (\mu^{(1)}, \ldots,\mu^{(s)})$.
Then we define the \textit{conjugate} of $ \Lambda$
via $\Lambda^\prime := (\blambda^\prime \mid \bmu^\prime ) $. 
Note that $ \Lambda^\prime \in {\cal L}_n(\alpha)$. 
Similarly, if $ \et = (\bT \mid \bu) \in  \Tab(\Lambda)$
we define the conjugate tableau $ \et^{\prime} \in \Tab(\Lambda^{\prime})$ 
by conjugating all the components of $ \et $, mimicking what we did for $ \Lambda $. 
We have the following compatibility between the dominance orders and conjugation: 
$ \Lambda \unlhd \overline{\Lambda} $ iff $ \overline{\Lambda}^\prime \unlhd {\Lambda}^\prime $ and
$ \es \unlhd \et $ iff $ \et^{\prime} \unlhd \es^{\prime} $.

\medskip

The following simple Lemma shall be used repeatedly. 
\begin{lem}{\label{simpleLemma}}
  Let $ \alpha \in \Par$. Suppose that $ A \in  \mathcal{SP}_n $ is of type $ \alpha$ 
  and that $ i $ and $i+1 $ are in distinct blocks
  of $ A$. Then in $ \Ea$ and $ \Eak $ the generator $g_i $ verifies the symmetric group quadratic relation
  when acting on $ { \be}_{A} $, that is 
  \begin{equation} { \be}_{A} g_i^2 = { \be}_{A}.
   \end{equation} 
\end{lem}
\begin{demo}
  By relation (\ref{E9}) we have that
\begin{equation}
{ \be}_{A} g_i^2 = { \be}_{A} (1+ (q-q^{-1})e_i g_i )=   \be_{A} + (q-q^{-1}) \be_{A}e_i g_i = \be_{A}
\end{equation}  
where we used (3) of Proposition \ref{commuting} for the last equality.
\end{demo}

\medskip
For $ \Lambda = (\blambda \mid \bmu ) \in {\cal L}_n(\alpha)$, and $ \Lambda$-tableaux
$\es = (\Bs \mid \bu) $ and $ \et = (\bT \mid \bv) $, we define elements of $ \Ea $ or $ \Eak$ 
\begin{equation}\label{multiMMurphy}
  x_{\Bs \bT}:=  g_{d(\Bs)}^{\ast} { \be}^{}_{\blambda} x^{}_{\blambda} g^{}_{d(\bT)},  \, \,\, \,\, \,\, \,\, \,\, \,
  y_{\Bs \bT}:=  g_{d(\Bs)}^{\ast} { \be}^{}_{\blambda}y^{}_{\blambda} g^{}_{d(\bT)}.
\end{equation}

We also need the following Lemma. 
\begin{lem}{\label{firstduality}}
  Suppose that $ \Lambda =(\blambda \mid \bmu )   \in {\cal L}_n(\alpha)  $, that $  \Lambda_1 =
(\blambda_1 \mid \bmu_1 ) 
  \in {\cal L}_n(\alpha)$ and 
that $ \Bs, \bT \in \std(\blambda) $ and $ \Bs_1, \bT_1 \in \std(\blambda_1) $.
Suppose moreover that 
$  x_{\Bs\bT} y_{ \Bs_1 \bT_1     } \neq 0 $.  Then $ \bT \unlhd \Bs_1^{\prime}          $.
\end{lem}
\begin{demo}
  For usual integer partitions this result is well-known, and as we shall see, the proof can be
  reduced to that case. 
Without loss of generality we may assume that $ \Bs = \bT^{\blambda} $ and
$ \bT_1 = \bT^{\blambda_1} $ and so we have that 
\begin{equation}\label{47}
\begin{array}{lr}
 x^{}_{\Bs\bT} y^{}_{ \Bs_1 \bT_1     } = 
 x_{\bT^{\blambda}  \bT} y_{  \Bs_1   \bT^{\blambda_1} } = \be^{}_{\blambda} x^{}_{\blambda} g^{}_{d(\bT)}
 g_{d( \Bs_1)}^{\ast}  \be^{}_{\blambda_1}   y^{}_{\blambda_1}
=  x^{}_{\blambda} g^{}_{d(\bT)} 
g_{d( \Bs_1)}^{\ast}    \be^{}_{A_{\blambda} d(\bT)^{} d( \Bs_1)^{-1}  }  \be^{}_{{\blambda}}
y^{}_{\blambda_1}. 
\end{array}
\end{equation}
For the last step
we used that $ x^{}_{\blambda} \be^{}_{{\blambda}} = \be^{}_{{\blambda}}  x^{}_{\blambda} $,
together with part (2) of Proposition \ref{commuting} 
and the fact that $ \be^{}_{{\blambda}} = \be^{}_{{\blambda_1}} $ since
$ \blambda $ and $ \blambda_1 $ are both of type $ \alpha$.
Now by hypothesis,  (\ref{47}) is nonzero and so 
$   \be_{A_{\blambda} d(\bT)^{} d( {\Bs_1})^{-1}  }  \be_{\blambda} \neq 0 $.  Hence using
the orthogonality of the $ \be_A$'s we get that 
\begin{equation}\label{andsoon} A_{\blambda} d(\bT) = A_{\blambda} d({\Bs}_1).
\end{equation}
  
We now choose decompositions with respect to $ \Si_{  \norm{\blambda} } = \Si_{\alpha} $
as in (\ref{decompositionintialkind})
\begin{equation}
 d(\bT ) = d(\bT_0) w_{\bT} \, \, \, \, \, \, \mbox{and} \, \, \, \, \, \,
 d(\Bs_1 ) = d((\Bs_1)_0) w_{\Bs_1}.
\end{equation}
Clearly $ \Si_n $ acts transitively on $ \mathcal{SP}_n^{\alpha} $
with $ \Si_{\alpha} $ as the stabilizer group of $ A_{\blambda}$ and so we have an identification
$ \mathcal{SP}_n^{\alpha} = \Si_{\alpha} \backslash \Si_n $. 
From (\ref{andsoon}) we have that 
$ \Si_{\alpha} w_{\bT} = \Si_{\alpha} w_{\Bs_1} $ and so 
$  w_{\bT} = w_{\Bs_1} $.
Hence (\ref{47}) becomes
\begin{equation}\label{49}
 \be^{}_{\blambda} x^{}_{\blambda} g^{}_{d(\bT_0)} g^{}_{  w_{\bT} } g^{\ast}_{  w_{\bT} } 
 g_{d( (\bT_1)_0)}^{\ast}     y^{}_{\blambda_1} = 
  x^{}_{\bT^{\blambda}  \bT_0} \be^{}_{\blambda} g^{}_{  w_{\bT} } g^{\ast}_{  w_{\bT} } 
y^{}_{  (\Bs_1)_0   \bT^{\blambda_1 }}.
 \end{equation}
We can now use Lemma {\ref{simpleLemma}} repeatedly on the simple transpositions that appear in a 
reduced expression for 
$  w_{\bT} $ and hence (\ref{49}) becomes
\begin{equation}
 \be^{}_{\blambda}
 x^{}_{\bT^{\blambda}  \bT_0} y^{}_{  (\Bs_1)_0   \bT^{\blambda_1 }}.
\end{equation}
In view of (\ref{thenweave}) and the fact that 
$ \bT_0 $ and $ (\Bs_1)_0 $ are both of the initial kind, 
the Lemma now follows from the similar result in the Hecke algebra setting.
\end{demo}  

\medskip
We wish to generalize the Lemma to a statement involving
$ m_{\es\et} $ and $ n_{\es\et}$. For this we need to recall from \cite{ER} some of the commutation relations
between the various factors of $ m_{\es\et} $ and $ n_{\es\et}$, as defined in
(\ref{maindefM}) and (\ref{maindefN}).

The commutation relations involving $ { \be}^{}_{\blambda} $ and
$ x_{\blambda} $ (resp. $ y_{\blambda} $)
are easy to describe,
since $ { \be}^{}_{\blambda} $ and
$ x_{\blambda} $ (resp. $ y_{\blambda} $)
commute with all the factors of 
(\ref{maindefM}) (resp. (\ref{maindefN})), except the ones on the extreme left and right.
The commutation relations between
$ \BB^{}_{d(\bv)} $ and $ g^{}_{d(\bT)} $ are more complicated to describe, but we shall here only
need the case when $ \bT $ is of the initial kind, that is $ \bT= \bT_0$. 
Combining Lemma 21 and Lemma 22 of
\cite{ER} we get in that case that
\begin{equation}\label{morecommutation}
{ \be}^{}_{\blambda}  \BB^{}_{d(\bv)}  g^{}_{d(\bT_0)}=
{ \be}^{}_{\blambda}  g^{}_{d(\bT_0)} \BB^{}_{d(\bv)}  
\end{equation}
and similarly for the two factors on the left of 
(\ref{maindefM}) and (\ref{maindefN}). 
We can now prove the following statement. 
\begin{lem}{\label{secondduality}}
  Suppose that $ \Lambda =(\blambda \mid \bmu )  $ and that $  \Lambda_1 =
(\blambda_1 \mid \bmu_1 ) 
  \in {\cal L}_n(\alpha)$ and 
that $ \es, \et \in \std(\Lambda) $ and $ \es_1, \et_1 \in \std(\Lambda_1) $.
If in $ \Ea $ or $ \Eak $ we have 
$  m_{\es \et } n_{ \es _1 \et_1     } \neq 0 $ then $ \et \unlhd \es_1^{\prime}          $.
\end{lem}
\begin{demo}
Let $\es = (\Bs \mid \bu) $, $ \et = (\bT \mid \bv) $, 
$\es_1 = (\Bs_1 \mid \bu_1) $ and $ \et_1 = (\bT_1 \mid \bv_1) $.
As in the proof of Lemma {\ref{firstduality}}, we can 
without loss of generality assume that $ \es = \et^{\Lambda} $ and
$ \et_1 = \et^{\Lambda_1} $ and so we have that 
\begin{equation}\label{product}
  m_{\es \et } n_{ \es _1 \et_1     } =
 { \be}^{}_{\blambda} x^{}_{\blambda} b^{}_{\bmu}\BB^{}_{d(\bv)}  \,
g^{}_{d(\bT)}
g_{d(\Bs_1)}^{\ast}  
\BB_{d(\bu_1)}^{\ast}  c^{}_{\bmu_1} y^{}_{\blambda_1} { \be}^{}_{\blambda_1}. 
\end{equation}
We now argue largely as in Lemma {\ref{firstduality}}.  
We first observe that $ { \be}^{}_{\blambda_1} = { \be}^{}_{\blambda} $
since both $ \blambda$ and $ \blambda_1 $ are of type $ \alpha$.
Using the commutation rules
involving $ { \be}^{}_{\blambda} $, see the paragraph prior to (\ref{morecommutation})
and part (3) of Proposition \ref{commuting}, 
we get that (\ref{product}) is equal to
\begin{equation}\label{wemauas}
 x^{}_{\blambda} b^{}_{\bmu}\BB^{}_{d(\bv)}  \,
g^{}_{d(\bT)}
g_{d(\Bs_1)}^{\ast}  
\BB_{d(\bu_1)}^{\ast}  
{ \be}^{}_{A_{\blambda} d(\bT) d(\Bs_1)^{-1} } { \be}^{}_{\blambda} c^{}_{\bmu_1}
y^{}_{\blambda_1}. 
\end{equation}
By hypothesis, (\ref{product}) is nonzero and so 
$ { \be}^{}_{A_{\blambda} d(\bT) d(\Bs_1)^{-1} }=  { \be}^{}_{\blambda} $ 
which implies that 
$ A_{\blambda} d(\bT) = A_{\blambda} d({\Bs}_1) $.
As in Lemma {\ref{firstduality}} we have decompositions
$ d(\bT ) = d(\bT_0) w_{\bT} $ and
$ d(\Bs_1 ) = d((\Bs_1)_0) w_{\Bs_1} $ 
and so 
$ A_{\blambda} w_{\bT} = A_{\blambda} w_{\Bs_1} $ and
$  w_{\bT} = w_{\Bs_1} $.
Hence (\ref{wemauas}) becomes
\begin{equation}\label{68becomes}
\begin{split}
{ \be}^{}_{\blambda}  x^{}_{\blambda} b^{}_{\bmu}\BB^{}_{d(\bv)}  \,
g^{}_{d(\bT_0)} g^{}_{  w_{\bT} } g^{\ast}_{  w_{\bT} }  
g_{d((\Bs_1)_0)}^{\ast}  
\BB_{d(\bu_1)}^{\ast} c^{}_{\bmu_1}  
y^{}_{\blambda_1}= \\
{ \be}^{}_{\blambda}   x^{}_{\bT^{\blambda}  \bT_0} b^{}_{\bmu}\BB^{}_{d(\bv)}  \,
 g^{}_{  w_{\bT} } g^{\ast}_{  w_{\bT} }  
\BB_{d(\bu_1)}^{\ast} c^{}_{\bmu_1}  
y^{}_{  (\Bs_1)_0   \bT^{\blambda_1 }} =\\
  x^{}_{\bT^{\blambda}  \bT_0} b^{}_{\bmu}\BB^{}_{d(\bv)}  \,
{ \be}^{}_{\blambda } g^{}_{  w_{\bT} } g^{\ast}_{  w_{\bT} }  
\BB_{d(\bu_1)}^{\ast} c^{}_{\bmu_1}  
y^{}_{  (\Bs_1)_0   \bT^{\blambda_1 }}.
\end{split}
\end{equation}
Using Lemma {\ref{simpleLemma}} repeatedly on a reduced expression for $  w_{\bT}  $
we may cancel $  g^{}_{  w_{\bT} } g^{\ast}_{  w_{\bT} }   $ out and so (\ref{68becomes})
becomes 
\begin{equation}\label{53}
{ \be}^{}_{\blambda}     x^{}_{\bT^{\blambda}  \bT_0} b^{}_{\bmu}\BB^{}_{d(\bv)}  \,
\BB_{d(\bu_1)}^{\ast} c^{}_{\bmu_1}  
y^{}_{  (\Bs_1)_0   \bT^{\blambda_1 }} =
{ \be}^{}_{\blambda}     x_{\bT^{\blambda}  \bT_0} y_{  (\Bs_1)_0   \bT^{\blambda_1 }}
b^{}_{\bmu}\BB^{}_{d(\bv)}  \,
\BB_{d(\bu_1)}^{\ast} c^{}_{\bmu_1}
\end{equation}
where we used that $  y^{}_{  (\Bs_1)_0   \bT^{\blambda_1 }} $ is of the initial kind. 
We now argue once again as in Lemma {\ref{firstduality}}.
We know that (\ref{53}) is nonzero. 
Since $ \bT_0 $ and $ (\Bs_1)_0 $ are of the initial kind we 
deduce from Lemma \ref{firstduality}
that $ \bT \unlhd \Bs_1^{\prime} $. But since (\ref{53}) is nonzero we also get that 
$   b^{}_{\bmu}\BB^{}_{d(\bv)}  \,
\BB_{d(\bu_1)}^{\ast}  
 c^{}_{\bmu_1} $ is nonzero. Moreover
$ \bv $ and $ \bu_1 $ are of the initial kind as well, this time by construction, and so we conclude that 
$ \bv \unlhd \bu_1^{\prime} $. This proves the Lemma.
\end{demo}

\section{Permutation modules} 
Suppose that $ \blambda \in \MC$ is of type $ \alpha$. Then 
we define the $\Ea$-\textit{permutation module }$ M(\blambda) $ as the following right
ideal of $\Ea$
\begin{equation}
M(\blambda) :=   \be_{\blambda} x_{\blambda} \Ea  \subseteq \Ea.
\end{equation}
We also have a specialized version of the permutation module that we denote the same way
\begin{equation}
M(\blambda) :=   \be_{\blambda} x_{\blambda} \Eak  \subseteq \Eak.
\end{equation}
$ M(\blambda) $ 
is a generalization of the permutation modules that appear in the representation theory
of the symmetric group, the Hecke algebra, the cyclotomic Hecke algebra, and so on, see 
\cite{DJM} and the references therein.
For any $ \Bs \in \rstd(\blambda)$, we introduce
$ x_{\Bs }:=   { \be}^{}_{\blambda} x^{}_{\blambda} g^{}_{d(\Bs)} \in M(\blambda)$
and define $ r_i =r_{i}^{\Bs}, c_i = c_i ^{{\Bs}},p_i= p_i^{\Bs} $ via
$ \Bs(r_i, c_i, p_i) =  i $.
We then have the following Lemma which describes the $ \Ea$-module (or $ \Eak$-module) structure on 
$ M(\blambda)  $.

\begin{lem}\label{descr}
  \begin{itemize}
\setlength\itemsep{-1.5em}
\item[(1)] The set $\{ x_{\Bs } \mid \Bs \in \rstd(\blambda) \} $ is a basis for $ M(\blambda)  $,
  over $ S $ or $ \kk$.
  \newline
\item[(2)]
  Let $  \Bs  \in \rstd(\blambda)$. Then
\begin{equation}{\label{actiongi}}
\begin{array}{l}
  x_{\Bs}  g_i=\left\{\begin{array}{ll} x_{\Bs s_i}    &\mbox{ if } p_i \neq p_{i+1}  \\
q x_{\Bs }         &\mbox{ if } p_i = p_{i+1} \mbox{ and } r_i= r_{i+1}  \\
 x_{\Bs s_i}     &\mbox{ if }p_i = p_{i+1} \mbox{ and } r_i < r_{i+1} \\
(q-q^{-1}) x_{\Bs }     +  x_{\Bs s_i}     &\mbox{ if }
p_i = p_{i+1} \mbox{ and } r_i > r_{i+1}
\end{array}\right.
\\ \\
x_{\Bs}  e_i=\left\{\begin{array}{ll} x_{\Bs}    &\mbox{ if } p_i = p_{i+1}  \\
0         &\mbox{ otherwise.} 
\end{array}\right.
\end{array}
\end{equation}
\end{itemize}
\end{lem}
\begin{demo}
Using (\ref{basessubalgebras}), we get that 
the set $ \{ \be_{\blambda} g_w \, | \, w \in \Si_n \}$ generates $ M(\blambda) $ over $ S$. 
On the other hand, for $ w \in \Si_n $ we have the decomposition 
$ w = w_0 d(\Bs) $ with $ w_0 \in \Si_{\blambda} $ and $ \Bs \in \rstd(\blambda) $, and from this it follows 
via ({\ref{justasinYH}}) that 
\begin{equation}\label{66}
  \be_{\blambda}   x_{\blambda} g_w = q^{\ell(w_0)} \be_{\blambda} x_{\blambda} g_{\Bs} =
  q^{\ell(w_0)} x_{\Bs}.
\end{equation}  
Hence in fact $ \{x_{\Bs}\, |\,  \Bs \in \rstd(\blambda) \}$ generates $ M(\blambda)$ over
$ S$. 
In the expansion of $ x_{\Bs} $ in terms of $ \be_{\blambda} g_{w}^{\,\,\,\,\,\, ,} $s we know that
$ \be_{\blambda}g_{d(\Bs)} $
appears exactly once,
since the 
$ d(\Bs) $'s are distinguished coset representatives for $ \Si_{\blambda}$ in $ \Si_n$,
whereas the other terms are of the form $  \be_{\blambda} g_{d(\Bs)} g_{w} $ for $ w \in \Si_{\blambda}$. 
These elements all belong to the basis for $ \Ea $ given in (\ref{basessubalgebras}), and are all distinct, 
and so the $ x_{\Bs}$'s are linearly independent, proving (1).

We next prove (2), where we 
first consider $ x_{\Bs} e_i $. Here we have that 
\begin{equation} x_{\Bs} e_i = 
  \be_{\blambda} g_{d(\Bs)} e_i = g_{d(\Bs)} \be_{A_{\blambda} d(\Bs)}  e_i
\end{equation}  
and so the formula for $ x_{\Bs} e_i  $ follows from Proposition \ref{commuting}. 
Let us then  consider $ x_{\Bs} g_i $ when $ p_i = p_{i+1}$.
If $ r_i = r_{i+1} $ then we have $ d(\Bs) s_i = s_j d(\Bs) $
for some $ s_j \in \Si_{\blambda} $, 
and so $ x_{\Bs} g_i = q x_{\Bs}$
via (\ref{66}). If $ r_i < r_{i+1} $ then $ d(\Bs) s_i > d(\Bs)  $ and so
$ x_{\Bs} g_i =  x_{\Bs s_i}$
from the definition of $ x_{\Bs s_i}$. If $ r_i > r_{i+1} $
then $ d(\Bs) s_i < d(\Bs)  $. We can then choose a reduced expression for $ d(\Bs) $ ending
in $ s_i $ and so 
\begin{equation}\label{68}
  x_{\Bs} g_i = \be_{\blambda} g_{d(\Bs) } g_i = \be_{\blambda} g_{d(\Bs s_i) } g_i^2 =
  \be_{\blambda} g_{d(\Bs s_i) } \left(1 + (q-q^{-1}) e_i g_i \right) =    x_{\Bs s_i} + (q-q^{-1}) x_{\Bs }    
\end{equation}
where we used the quadratic relation (\ref{E9}) for the third equality, and the formula for $ x_{\bs} e_i $
for the last equality. Finally, if $ p_i \neq  p_{i+1}$ then either 
$ d(\Bs) s_i > d(\Bs)  $ or   $ d(\Bs) s_i < d(\Bs)  $ and so we can repeat the previous arguments, with
the only difference that the second term in (\ref{68}) disappears.
These arguments work over $ \kk $ as well and so the Lemma is proved,
\end{demo}
\medskip

There is another kind of permutation module
that we shall need.
Suppose that
$ \Lambda = (\blambda \mid  \bmu ) \in {\cal L}_n(\alpha) $. Then we define
the permutation module $ M(\Lambda) $ as the following right ideal of $ \Ea$
\begin{equation}
M(\Lambda) := \be_{\blambda} x_{\blambda} b_{\bmu} \Ea  \subseteq \Ea.
\end{equation}
Once again there is a specialized version of $ M(\Lambda) $ obtained by replacing
$ \Ea $ with $ \Eak$.
The only difference between $ M(\Lambda)  $ and $ M(\blambda) $ is the factor $ b_{\bmu} $
and so we have $ M(\Lambda)  \subseteq M(\blambda) $, with the inclusion being strict in general.
For $ M(\Lambda) $ we have the following Lemma. 

\begin{lem}\label{descr2}
Let $ \Lambda \in  {\cal L}_n(\alpha)$. 
Let $ m_{\et }$  be as in (\ref{forsimplicity}). Then 
the set $\left\{ m_{\et } \mid \et \in \rstd(\Lambda) \right\} $ is an $S$-basis (or a $ \kk $-basis)
for $ M(\Lambda)  $.
\end{lem}
\begin{demo}
 For $ w \in \Si_n $ we have the decomposition
 $ w = w_0 d(\bT) $ with $ w_0 \in \Si_{\blambda}$ and $ \bT \in \rstd(\blambda) $ and hence,
 arguing as in (\ref{66}), we get 
\begin{equation}\label{77}
  \be_{\blambda} x_{\blambda} b_{\bmu} g_w =
  \be_{\blambda} x_{\blambda} b_{\bmu} g_{w_0} g_{\bT} =
  \be_{\blambda} x_{\blambda}g_{w_0} b_{\bmu}  g_{\bT} =
q^{\ell(w_0)}\be_{\blambda} x_{\blambda} b_{\bmu}  g_{\bT}.
\end{equation}  
Hence
$ \left\{  \be_{\blambda}x_{\blambda} b^{}_{\bmu}  g_{\bT}
\, \big| \, \bT \in \rstd(\blambda) \right\}$ generates $ M(\Lambda) $ over $S$, but 
the appearing $ \bT $ may not be increasing 
in the sense of the definition of $ \rstd(\Lambda) $.  On the other hand, if $ \bT $
is not increasing we can find 
$ B_u \in  \Si_m$
such that the multitableau
$ u \ast \bT  \in  \rstd(\Lambda) $ defined 
by
\begin{equation}\label{ast defi}
  d( u \ast \bT ) =    B_u d(\bT)
\end{equation}
becomes increasing; here $ \Si^m_{\Lambda}
= \Si_{m_1} \times \ldots \times \Si_{m_s} $ is the subgroup of 
$ \Si_n$ 
 introduced in (\ref{abstractgroup}). We then get
\begin{equation}\label{78}
  \be_{\blambda} x_{\blambda} b_{\bmu}  g_{\bT} =
  \be_{\blambda} x_{\blambda} b_{\bmu} \BB_u^{-1} \BB_u   g_{d(\bT)} =  
  \be_{\blambda} x_{\blambda} b_{\bmu} \BB_u^{-1}  g_{  B_u   d(\bT) } =
  \be_{\blambda} x_{\blambda} b_{\bmu} \BB_u^{-1} g_{   u \ast \bT }    
\end{equation}  
where we for the second equality used Lemma 54 of \cite{ER}.
Suppose $ \bmu = (\mu^{(1)}, \ldots, \mu^{(s)}) $ and denote by 
$ \Si^m_{\Lambda, \bmu} :=  \Si_{ \mu^{(1)} }\times \ldots \times \Si_{ \mu^{(s)} }$
the corresponding diagonal Young subgroup of  $ \Si^m_{\Lambda}$. It gives 
rise to a diagonal decomposition of $ B_u^{-1} $ that is $B_u^{-1}= (B_u)_0 d(\bu) $
where $ (B_u)_0 \in \Si^m_{\Lambda, \bmu} $
and $ \bu $ is a row standard $\bmu $-multitableau of the initial kind. 
Setting $ \et = ( u \ast \bT \mid \bu ) \in {\cal L}_n(\alpha)$, we can now rewrite (\ref{78}) as
\begin{equation}\label{79} 
  \be_{\blambda} x_{\blambda} b_{\bmu} \BB_u g_{   u \ast \bT }   =
  \be_{\blambda} x_{\blambda} b_{\bmu} \BB_{\bu} g_{   u \ast \bT }    = m_{\et}.
\end{equation}   
Thus we have proved that the set $ \{m_{\et}\}$ generates $ M(\Lambda) $.
To show that it is linearly independent we argue as follows.
Suppose that $ \es = ( \Bs \mid \bu) $, that is 
\begin{equation}
   m_{\es}=   { \be}^{}_{\blambda} x^{}_{\blambda}  b^{}_{\bmu}\BB^{}_{d(\bu)}  \,
g^{}_{d(\Bs)}.
\end{equation}
We know that $ \bu$ is row standard of the initial kind and therefore
$  b^{}_{\bmu}\BB^{}_{d(\bu)}  $ expands as a sum of terms of the form $ \BB^{}_{u_0d(\bu)} $ 
where $ u_0 \in  \Si^m_{\Lambda, \bmu} $. But using Lemma 54 of \cite{ER} we get that 
\begin{equation} \BB^{}_{u_0d(\bu)} g^{}_{d(\Bs)} =
 g^{}_{B_{u_0} d(\bu)    d(\Bs)}.
\end{equation}
We now note that the multitableau $ \Bv $ corresponding to $ B_{u_0} d(\bu)    d(\Bs) $, in other words
$ \Bv $ given by
$ d(\Bv) = B_{u_0} d(\bu)    d(\Bs) $, is a row standard $\blambda$-tableau since it is obtained from
$ \bT $ by permuting components with equal $ \lambda^{(i)}$'s and 
so $ m_{\es} $ expands as a sum of terms of the
form $ { \be}^{}_{\blambda}   g^{}_{x B_{u_0} d(\bu)    d(\Bs)} $ where
$ u_0 \in  \Si^m_{\Lambda, \bmu} $ and 
$ x \in \Si_{\blambda} $. These expansions are distinct for distinct pairs $  ( \Bs \mid \bu) $
and so the set $ \{m_{\es}\}$ is linearly independent, as claimed. Hence it is an $S$-basis,
and therefore also a $ \kk $-basis, for $M(\Lambda) $.

\end{demo}

\medskip
We next give a description of the action of $ \Ea $ on $ \{m_{\es}\}$. 
We first introduce the following useful notation. Let $ \Bs \in  \rstd(\blambda) $ and let 
let $ s_i \in \Sigma_n \subseteq \Si_n$ be a simple transposition. Then we define $  \Bs \cdot s_i $ via
\begin{equation}
  \Bs \cdot s_i := 
\left\{\begin{array}{ll}
\Bs  & \mbox{if  } r_i = r_{i+1}  \\
\Bs s_i & \mbox{if  } r_i \neq r_{i+1}.
  \end{array}
  \right.
\end{equation}
This extends to an action of $ \Si_n $ on $ \rstd(\blambda) $ that we denote $ (\Bs, w ) \mapsto \Bs \cdot w $.
We extend this further to $ {\cal L}_n(\alpha) $ as follows. Let $ \es =
(\Bs \mid \bu) \in  {\cal L}_n(\alpha) $ and suppose that 
$ Shape(\T^{(p_i)}) = Shape(\T^{(p_{i+1})}) $. Then we define a transposition 
$  \tau_i \in \Si_{\Lambda}^m$ via
$ \tau_i = (p_i, p_{i+1})$. 
For $ \Bs \in \rstd(\blambda) $ we set
$ m_i :=  \min(Shape(\T^{(p_i)}) $.
Using the notation $ \tau_i \ast \Bs  $ introduced in
(\ref{ast defi}) we then define

\begin{equation}\label{extendstoan}
  \es \cdot  s_i := \left\{
      \begin{array}{ll}
(\Bs \cdot s_i\mid \bu) & \mbox{if  } p_i = p_{i+1}   \mbox{ or otherwise    }  \\ (\Bs \cdot s_i\mid \bu)
& \mbox{if } Shape(\T^{(p_i)}) \neq Shape(\T^{(p_{i+1})})   \\ (\Bs \cdot s_i\mid \bu)
 & \mbox{if }  Shape(\T^{(p_i)}) = Shape(\T^{(p_{i+1})})
  \mbox{ and  }  \{ m_i , m_{i+1} \} \neq \{i,i+1\}  \\ 
  (\tau_i \ast(\Bs \cdot \tau_i) \mid \bu \cdot \tau_i)  &
  \mbox{if }   Shape(\T^{(p_i)}) = Shape(\T^{(p_{i+1})})
  \mbox{ and  }  \{ m_i , m_{i+1} \} = \{i,i+1\} .
  \end{array}
  \right.
\end{equation}

One checks that these formulas
extend to an action of $ \Si_n $ on $ \rstd(\Lambda) $.
They have their origin in the straightening procedure
that carries $ \Bs s_i $ to an increasing row standard tableaux. Here are two examples,
corresponding to the last case in (\ref{extendstoan}):\begin{equation}\begin{array}{c}
{\footnotesize\left(\;\left(\;\young(15)\;,\young(26)\;,\young(34)\;\right) \big|
  \left(  \young(12,3) \right)\;\right)} \cdot s_1 =
{\footnotesize\left(\;\left(\;\young(16)\;,\young(25)\;,\young(34)\;\right) \big|
  \left(  \young(12,3) \right)\;\right)}
\\[4mm]
{\footnotesize\left(\;\left(\;\young(15)\;,\young(26)\;,\young(34)\;\right) \big|
  \left(  \young(13,2) \right)\;\right)} \cdot s_2 =
{\footnotesize\left(\;\left(\;\young(15)\;,\young(24)\;,\young(36)\;\right) \big|
  \left(  \young(23,1) \right)\;\right)}.
\end{array}\end{equation}

\begin{lem}\label{descr2A}
  Suppose that $  \es = (\Bs \mid \bu ) \in \rstd(\Lambda)$. Then 
  \begin{equation}{\label{actiongiaa}}
\begin{array}{l}    
  (1) \, \, \,  m_{\es}   g_i=\left\{
\begin{array}{ll} q m_{\es \cdot s_i}    & \mbox{ if } p_i = p_{i+1} \mbox{ and } r_i = r_{i+1}   \\
  m_{\es \cdot s_i}    & \mbox{ if } p_i = p_{i+1} \mbox{ and } r_i < r_{i+1}   \\
  (q-q^{-1}) m_{\es } +  m_{\es \cdot s_i}    & \mbox{ if } p_i = p_{i+1} \mbox{ and } r_i > r_{i+1}   \\
m_{\es \cdot s_i}    & \mbox{ if } p_i \neq p_{i+1}   \\    
\end{array}\right.
  \\ \\
(2)  \, \, \,   m_{\es}   e_i=\left\{\begin{array}{ll}   m_{\es}     &\mbox{ if } p_i = p_{i+1}  \\
0        &\mbox{ if }  p_i \neq p_{i+1}.
  \end{array}\right. \end{array}
\end{equation}
  \end{lem}
\begin{demo}
  The proof of (2) is analogous to the corresponding statement for
  $ x_{\Bs } e_i $ in Lemma \ref{descr}, in detail 
\begin{equation}  
  m_{\es} e_i = \be_{\blambda} x_{\blambda} b_{\bmu} \BB_{\bu} g_{d(\Bs)}  e_i = 
 x_{\blambda} b_{\bmu} \BB_{\bu} g_{   d(\Bs) }
\be_{A_{\blambda} d(\Bs)} e_i 
\end{equation}
and we conclude via Proposition \ref{commuting}.

We therefore consider (1) where we first focus on the case $ p_i = p_{i+1}$. 
If $ r_i = r_{i+1} $ we have $ d(\Bs) s_i = s_j d(\Bs) $
for some $ s_j \in \Si_{\blambda} $, and so 
\begin{equation}\label{85}  
  m_{\es} g_i = \be_{\blambda} x_{\blambda} b_{\bmu} \BB_{\bu} g_{d(\Bs)}  g_i =  \be_{\blambda}
 x_{\blambda} b_{\bmu} \BB_{\bu} g_j g_{   d(\Bs) }.
\end{equation}
But $ s_j $ is of the form $s_j = d(\bT) $ where  $\bT $ is a multitableau of the initial kind and so (\ref{85})
becomes
\begin{equation}\label{86}  
  x_{\blambda} g_j  b_{\bmu} \BB_{\bu} g_{   d(\Bs) } = q  x_{\blambda}   b_{\bmu} \BB_{\bu} g_{   d(\Bs) } = 
  q m_{ \es},  
\end{equation}
as claimed. The two other cases when $ p_i = p_{i+1} $, that is $ r_i > r_{i+1} $ or $ r_i < r_{i+1} $, are
proved the same way as in Lemma \ref{descr}. 

We then finally consider the case $ p_i \neq p_{i+1}$. Arguing as in Lemma \ref{descr} we have
here 
\begin{equation}\label{87}  
  m_{\es} g_i = \be_{\blambda} x_{\blambda} b_{\bmu} \BB_{\bu} g_{d(\Bs)}  g_i =
\be_{\blambda} x_{\blambda} b_{\bmu} \BB_{\bu} g_{d(\Bs  s_i)}   
  \end{equation}
which is equal to $ m_{\es \cdot s_i } $ in all cases except $ \{ m_i, m_{i+1} \} = \{ i, i+1\} $
where $ \bT $ is not increasing. But in that case, arguing as we did 
for (\ref{78}) and (\ref{79})
we can rewrite (\ref{87})
as follows
\begin{equation}
  \be_{\blambda} x_{\blambda} b_{\bmu} \BB_{\bu} g_{d(\Bs  s_i)} =
  \be_{\blambda} x_{\blambda} b_{\bmu} \BB_{\bu} \BB_{\tau_i} \BB_{\tau_i} g_{d(\Bs  s_i)} = 
  \be_{\blambda} x_{\blambda} b_{\bmu} \BB_{\bu  \cdot \tau_i }   g_{d( \tau_i \ast (\Bs  s_i))} = 
m_{\es \cdot s_i }
\end{equation}
as claimed.
\end{demo}  

\medskip

Given the basis for $ M(\Lambda) $, we now introduce a bilinear form $ (\cdot, \cdot)_{\Lambda} $
on $ M(\Lambda) $ as follows
\begin{equation}\label{bildef}
  (m_{\es}, m_{\et})_{\Lambda} := \left\{ \begin{array}{ll} 1 & \mbox{if } \es = \et\\
    0 & \mbox{otherwise. } \end{array} \right.
\end{equation}

It is a generalization of the bilinear forms on the classical permutation modules,
see for example \cite{DJM}.
\begin{lem}\label{lemma13}
The following statements hold. 
  \begin{itemize}
\setlength\itemsep{-1.5em}
\item[(1)] 
$ (\cdot, \cdot)_{\Lambda} $ is symmetric and nondegenerate.
  \newline
\item[(2)]
  $ (\cdot, \cdot)_{\Lambda} $ is invariant in the sense that for all $ m,m_1 \in M(\Lambda) $ and
  $ a \in \Ea $ or $ a\in \Eak $ we have that 
$ (m a, m_1)_{\Lambda} = (m , m_1 a^{\ast})_{\Lambda} $. 
\end{itemize}
\end{lem}

\begin{demo}
  (1) follows from the fact that $\{ m_{\es} \} $ is an orthogonal basis for $ M(\Lambda) $. To show (2) we may assume that $ m= m_{\es}, m_1 = m_{\et}   $ and $ a = g_i$ or $ a = e_i$
  and must check $ ( m_{\es} a,  m_{\et})_{\Lambda} = ( m_{\es} , m_{\et }a )_{\Lambda} $ for all
  possibilities of $ r_j^{\Bs}, c_j^{\Bs}, p_j^{\Bs}, r_j^{\bT}, c_j^{\bT}, p_j^{\bT} $ and $ j = i, i+1$.
Here the number of cases is reduced to the half by the symmetry of $ (\cdot, \cdot)_{\Lambda} $.
Moreover if $ p_i^{\Bs} = p_{i+1}^{\Bs} $ and
$ p_i^{\bT} \neq p_{i+1}^{\bT} $ then it follows immediately from Lemma \ref{descr2A} that
$  (m_{\es}e_i , m_{\et})_{\Lambda} = (m_{\es} , m_{\et} e_i)_{\Lambda} =0$ and that
$  (m_{\es}g_i , m_{\et})_{\Lambda} = (m_{\es} , m_{\et} g_i)_{\Lambda} =0$
and so we only need
consider the cases where 
$ p_i^{\Bs} = p_{i+1}^{\Bs} $ and $ p_i^{\bT} = p_{i+1}^{\bT} $ or 
$ p_i^{\Bs} \neq p_{i+1}^{\Bs} $ and $ p_i^{\bT} \neq  p_{i+1}^{\bT} $. 
In the rest of the proof we shall use the formulas of Lemma \ref{descr2A} repeatedly.

Let us first consider $ a = g_i $ and 
suppose that 
$ p_i^{\Bs} = p_{i+1}^{\Bs} $ and $ p_i^{\bT} = p_{i+1}^{\bT} $. 
If $ r_i^{\Bs} < p_{i+1}^{\Bs} $ and $ r_i^{\bT} < p_{i+1}^{\bT} $ 
we have $ (m_{\es} g_i, m_{\et})_{\Lambda}
  = (m_{\es \cdot s_i}, m_{\et })_{\Lambda} = \delta_{ \es \cdot s_i , \et} $ whereas
  $ (m_{\es}  , m_{\et} g_i)_{\Lambda}
  = (m_{\es},  m_{\et} g_i  )_{\Lambda} = \delta_{ \es ,\et \cdot s_i } $, where $ \delta$ is the Kronecker delta, 
and so $ (m_{\es} g_i, m_{\et})_{\Lambda} = (m_{\es},  m_{\et} g_i  )_{\Lambda} $.

If $ r_i^{\Bs} < p_{i+1}^{\Bs} $ and $ r_i^{\bT} > p_{i+1}^{\bT} $
then again
$ (m_{\es} g_i, m_{\et})_{\Lambda}
 = \delta_{ \es \cdot s_i , \et} $
whereas $ (  m_{\es} g_i, m_{\et}           )_{\Lambda} = ( m_{\es},  m_{\et}   g_i)_{\Lambda} = 
(m_{\es }, m_{\et \cdot s_i} + (q-q^{-1})m_{\et  })_{\Lambda}
= \delta_{ \es  ,  \et \cdot s_i   }    $, since $ \es \neq \et $, and so 
the two sides coincide.

If $ r_i^{\Bs} > p_{i+1}^{\Bs} $ and $ r_i^{\bT} > p_{i+1}^{\bT} $
we have 
$ (m_{\es} g_i , m_{\et} )_{\Lambda} =
(m_{\es \cdot s_i } + (q-q^{-1})m_{\es }, m_{\et})_{\Lambda} = \delta_{\es \cdot s_i, \et } + (q-q^{-1}) \delta_{\es, \et}$
and $ (m_{\es} , m_{\et} g_i)_{\Lambda} = 
(m_{\es }, m_{\et \cdot s_i} + (q-q^{-1})m_{\et  })_{\Lambda}
= \delta_{ \es  , \et \cdot s_i} + (q-q^{-1}) \delta_{\es, \et}   $. Once again, the two sides coincide.

We next consider $ a = g_i $ but 
suppose that 
$ p_i^{\Bs} \neq  p_{i+1}^{\Bs} $ and $ p_i^{\bT} \neq p_{i+1}^{\bT} $.
Here we have $ (m_{\es} g_i , m_{\et} )_{\Lambda} =  (m_{\es \cdot s_i } , m_{\et} )_{\Lambda} =
\delta_{\es \cdot s_i , \et} $ and similarly 
$ (m_{\es}  , m_{\et}g_i )_{\Lambda} =  (m_{\es } , m_{\et  \cdot s_i}   )_{\Lambda} =
\delta_{\es  , \et \cdot s_i} $ and so the two side also coincide in this case.

\medskip
Finally we consider the case $ a  =e_i $ and must check that
$ (m_{\es } e_i , m_{\et  })_{\Lambda} = (x_{\es } , m_{\et  } e_i )_{\Lambda}  $.
But if $ p_i^{\Bs} = p_{i+1}^{\Bs} $ and $ p_i^{\bT} = p_{i+1}^{\bT} $
we have $ m_{\es } e_i = m_{\es } $ and $ m_{\et } e_i = m_{\et } $, and
otherwise, if $ p_i^{\Bs} \neq  p_{i+1}^{\Bs} $ and $ p_i^{\bT} \neq p_{i+1}^{\bT} $, 
we have $ m_{\es } e_i = 0 $ and $ m_{\et } e_i = 0 $. 
Hence the two sides also coincide in this case. This finishes the proof of the Lemma. 
\end{demo}

\medskip
We have the following crucial Lemma. 
\begin{lem}\label{samespirit}
Suppose that $ \Lambda =(\blambda \mid \bmu )  $ and 
that $ \es, \et \in \std(\Lambda) $. Then we have 
\begin{equation}
  ( m_{ \et } n_{ \et^{\prime} \es^{\prime}     } , m_{\es})_{\Lambda} \neq 0 \mbox{  or equivalently  }
( m_{ \et }  , m_{\es}    n_{ \es^{\prime} \et^{\prime}     })_{\Lambda} \neq 0  .
\end{equation}  
In particular
$ m_{ \et } n_{ \et^{\prime} \es^{\prime}     } \neq 0$ and
$  m_{\es}    n_{ \es^{\prime} \et^{\prime}     } \neq 0$.
\end{lem}
\begin{demo}
Suppose that $ \es = (\Bs \mid \bu) $ and $ \et = (\bT \mid \bv) $. 
Then the same 
chain of equalities that took us from 
(\ref{product}) to (\ref{53}) in Lemma \ref{secondduality},
but postmultiplied with
$ \BB^{}_{d(\bu^{\prime})} g^{}_{d(\Bs^{\prime})}$, 
gives us 
\begin{equation}\label{byreducingto}
\begin{array}{l}
  m_{\et } n_{ \et^{\prime} \es^{\prime}     } = 
\be^{}_{\blambda}
 x^{}_{\bT^{\blambda}  \bT_0} y_{ \bT_0^{\prime}     \bT^{\blambda^{\prime} }}
  b^{}_{\bmu}\BB^{}_{d(\bv)}  \,
\BB_{d(\bv^{\prime})}^{\ast} 
c^{}_{\bmu^{\prime}}  \BB^{}_{d(\bu^{\prime})} g^{}_{d(\Bs^{\prime})} = \\ 
\be^{}_{\blambda}
 x^{}_{\blambda} g_{w_{\blambda}} y_{ \blambda^{\prime}}
  b^{}_{\bmu}\BB^{}_{ w_{\bmu}}  \,
c^{}_{\bmu^{\prime}}  \BB^{}_{d(\bu^{\prime})} g^{}_{d(\Bs^{\prime})}
\end{array}
\end{equation}
where we used (\ref{Murphyduality}) for the last equality: note that
$ \bT_0 $ and $ \bv $ are both of the initial kind. 
Using the invariance property of the previous Lemma \ref{lemma13} we get from this that
\begin{equation}\label{93}
  \begin{array}{l}
 \setlength\itemsep{0.5em}   
  ( m^{}_{\et } n^{}_{ \et^{\prime} \es^{\prime}} \, , \, m^{}_{\es})_{\Lambda} =
 ( m^{}_{\et } n^{}_{ \et^{\prime} \es^{\prime}}
 g^{\ast}_{d(\Bs)} \BB^{\ast}_{d(\bu)} \, , \,  m^{}_{\et^{\Lambda}})_{\Lambda} = \\
 ( \be^{}_{\blambda}
 x^{}_{\blambda} g_{w_{\blambda}} y_{ \blambda^{\prime}}
  b^{}_{\bmu}\BB^{}_{ w_{\bmu}}  \,
c^{}_{\bmu^{\prime}}  \BB^{}_{d(\bu^{\prime})} g^{}_{d(\Bs^{\prime})}
 g^{\ast}_{d(\Bs)} \BB^{\ast}_{d(\bu)} \, , \, m^{}_{\et^{\Lambda}})_{\Lambda} 
\end{array} 
\end{equation}
In view of (\ref{22}), and the remarks preceding (\ref{22}), 
we have decompositions $ d(\Bs^{\prime}) = d(\Bs^{\prime}_0) d(\bT) $ and
$ d(\Bs) = d(\Bs_0) d(\bT) $ and hence we can use 
Lemma {\ref{simpleLemma}} repeatedly on a reduced expression of $ d(\bT) $ 
to rewrite $  \be^{}_{\blambda} g^{}_{d(\Bs^{\prime})} g^{\ast}_{d(\Bs)}
=  \be^{}_{\blambda} g^{}_{d(\Bs_0^{\prime})} g^{\ast}_{d(\Bs_0)} = 
 \be^{}_{\blambda} g^{\ast}_{w_{\blambda} } $ and so 
 (\ref{93}) becomes
 \begin{equation}\label{94}
 \begin{array}{l}
 \setlength\itemsep{0.5em}   
   ( \be^{}_{\blambda}
 x^{}_{\blambda} g_{w_{\blambda}} y_{ \blambda^{\prime}}^{} g^{\ast}_{w_{\blambda}} 
  b^{}_{\bmu}\BB^{}_{ w_{\bmu}}  
c^{}_{\bmu^{\prime}}  \BB^{\ast}_{w_{\bmu}}     , m^{}_{\et^{\Lambda}})_{\Lambda} =\\
   ( \be^{}_{\blambda}
x^{}_{\blambda} g_{w_{\blambda}} y_{ \blambda^{\prime}}^{}
  b^{}_{\bmu}\BB^{}_{ w_{\bmu}}  
  c^{}_{\bmu^{\prime}} ,  
\be^{}_{\blambda} x^{}_{\blambda}  g^{}_{w_{\blambda}}  b^{}_{\bmu}      \BB^{}_{w_{\bmu}}     )_{\Lambda}
=
( m_{\et_{\Lambda}} y^{}_{ \blambda^{\prime}} c^{}_{\bmu^{\prime}} 
      , m^{}_{\et_{\Lambda}})_{\Lambda} 
\end{array}
 \end{equation}
Note that (\ref{94}) does not depend on $ \es $ and $ \et $, only on $ \Lambda$.
The expressions in (\ref{94}) only involve multitableaux of the initial
kind and so the Lemma is reduced to the corresponding Hecke algebra statement.
To be precise, for $ \lambda \in \Par $, each element $w$ of
the double coclass $ \Si_{\lambda}  w_{\lambda} \Si_{\lambda^{\prime}}  $ has a unique
decomposition $ w = w_1 w_{\lambda} w_2 $ where $ w_1 \in  \Si_{\lambda} $ and 
$ w_2 \in  \Si_{\lambda^{\prime}} $ and moreover $ \ell(w) = \ell(w_1) + \ell(w_{\lambda}) + \ell(w_2) $,
see Lemma 1.6 of \cite{DJ}. Applying this on $ \Si_{\blambda} $ and $ \Si_{\bmu} $ we conclude that
$ m_{\et_{\Lambda}} $ occurs exactly once 
in the expansion of $ m_{\et_{\Lambda}} y^{}_{ \blambda^{\prime}} c^{}_{\bmu^{\prime}} $,
corresponding to the '1'-term in both $ y^{}_{ \blambda^{\prime}} $ and $c^{}_{\bmu^{\prime}} $,
and so (\ref{94}) is equal to 1.
The Lemma is proved.


\end{demo}

\section{The tensor space $ V^{\otimes n} $ module for $ \E$}
In this section we realize the tensor space $ V^{\otimes n} $ module for $ \E$, introduced
in \cite{Ry}, as a sum of the
permutation modules $ M(\blambda)$. When the dimension of $ V $ is sufficiently large,
it is known that $ V^{\otimes n} $ is a faithful $\E$-module, but here we are interested in the general case where
$ V^{\otimes n} $ may not be faithful. Using the results from the previous section we determine the annihilator
ideal in $ \E $ of $ V^{\otimes n} $. In turns out that is has a nice description in terms of the dual cellular basis
$ \{ n_{\es\et} \} $. It is a main point of our constructions and proofs that they work for
arbitrary $ \kk$.

\medskip
Let $V$ be the free $S$-module with basis
\begin{equation}\label{95}
{\mathcal B}:=  \{ v_{i}^s  \, | \,  1 \le i \le N, 1 \le s \le r \}
\end{equation}
that is 
$ V $ is of dimension $ rN$.
We then define linear maps $\mathbf{E}, \mathbf{G} \in \End_S(V^{\otimes 2})  $ via
\begin{equation}{\label{operatorE}}
  (v_i^s\otimes v_j^t)\mathbf{E}:=\left\{\begin{array}{ll}
  v_j^s\otimes v_i^t&\mbox{ if }\;s =  t\\
0&\mbox{ if }\; s \neq t \\
\end{array}\right.
\end{equation}
and
\begin{equation}{\label{operatorG}}
(v_i^s\otimes v_j^t)\mathbf{G}:=\left\{\begin{array}{ll}v_j^t\otimes v_i^s&\mbox{ if }\;s\neq t\\
qv_i^s\otimes v_j^t&\mbox{ if }\;s=t,\;i=j\\
v_j^t\otimes v_i^s&\mbox{ if }\;s=t,\;i <j\\
(q-q^{-1})v_i^s\otimes v_j^t+v_j^t\otimes v_i^s&\mbox{ if
}\;s=t,\;i >j.\end{array}\right.
\end{equation}

We extend them to linear  maps $\mathbf{E}_i$ and $\mathbf{G}_i$ acting
in the tensor space $V^{\otimes n}$ by letting
$\mathbf{E}$ and $\mathbf{G}$
act in the $i$'th and $i+1$'st
factors. The following is Theorem 1 of \cite{Ry} 
\begin{teo}
  The rules $ e_i \mapsto \mathbf{E}_i $ and $ g_i \mapsto \mathbf{G}_i $ endow
  $V^{\otimes n}$ with the structure of
  an $ \E$-module.
\end{teo}
In the case of
$ { \mathcal E}^{\C}_n(q) $ 
it was proved in \cite{Ry} that the specialized tensor module $ V^{\C, \otimes n} $ is faithful when
$  r,N \ge n $. In \cite{ER} this faithfulness statement was generalized to $ V^{\otimes n} $ itself, 
but still only for $  r,N \ge n $; in fact it does not hold otherwise, as we shall see.

\medskip
Let $ \seqN $ be the set of sequences $ \underline{i}= (i_1, \ldots, i_n ) $ of
integers where each $ i_j $ belongs to $ \{ 1,\ldots, N\} $, and let similarly 
$ \seqr $ be the set of sequences $ \underline{s}= (s_1, \ldots, s_n ) $ where each 
$s_j $ belongs to $ \{ 1,\ldots, r\} $. For such $ \underline{i} $ and
$ \underline{s} $ we define
\begin{equation}
v_{\underline{i}}^{\underline{s}}:=v_{i_1}^{s_1} \otimes \cdots  \otimes  v_{i_n}^{s_n}\in V^{\otimes n}.
\end{equation}
Then the set $ \{ v_{\underline{i}}^{\underline{s}} \mid \underline{i} \in \seqN,  \underline{s} \in \seqr \} $
is a basis for $ V^{\otimes n} $.

\medskip
To $ \underline{s} \in \seqr $ we associate sets $ I_j :=\{ i \, | \, s_i =j \} $, for $ j=1,\ldots, r$. 
The $I_j$'s may be empty, but leaving out the empty $ I_j$'s we obtain a 
set partition in $ \mathcal{SP}_n $, that we denote $ A_{\underline{s}} $, and we say that $  \underline{s}  $ is of type $ \alpha $ if $ A_{\underline{s}} $
is of type $ \alpha$. For example, if $ r = 4,n=13 $ and $ \underline{s} = (1,2,2,1,2,2,2,1,2,4,1,2,2) $ then
$A_{\underline{s}} = \left\{ \{1,4,8,11\}, \{ 2, 3,5,6,7,9,12,13 \}, \{ 10 \} \right\} $, and so $ \underline{s}$ is
of type $ (8,4,1)$.

\medskip
Recall now the idempotent decompositions $ \sum_{ A \in \mathcal{SP}_n}  \mathbb{E}_A = 1 $
and $ \sum_{\alpha\in\Par }  \mathbb{E}_{\alpha} = 1 $.
They give rise to decompositions of $ V^{\otimes n} $
\begin{equation}\label{103}
  V^{\otimes n} = \bigoplus_{A \in \mathcal{SP}_n } \be_A   V^{\otimes n},\, \, \, \, \, \, \, \, \, \, \, \, \, \, \, 
V^{\otimes n} = \bigoplus_{\alpha\in\Par} \be_\alpha   V^{\otimes n}.
\end{equation}
Note that $ \be_\alpha   V^{\otimes n} $ is canonically an $ \Ea$-module.

\medskip
We have the following Lemma which gives precise descriptions of $ \be_A   V^{\otimes n} $
and $ \be_\alpha   V^{\otimes n}$.
\begin{lem}\label{16}
With the above notations the following statements hold 
  \begin{itemize}
\setlength\itemsep{-0.5em}
\item[(1)] $ E_A v_{\underline{i}}^{\underline{s}} = \left\{ \begin{array}{ll}
v_{\underline{i}}^{\underline{s}} & \mbox{ if } A \subseteq  A_{\underline{s}} \\
0 & \mbox{ otherwise }  \end{array} \right.$
  \newline
\item[(2)] $ \be_A   V^{\otimes n} = {\rm span}_S
\{  v_{\underline{i}}^{\underline{s}} \mid   A_{\underline{s}} = A  \} $
  \newline
\item[(3)] $  \be_{\alpha}   V^{\otimes n} = {\rm span}_S
\{  v_{\underline{i}}^{\underline{s}} \mid  \underline{s} \mbox{ is of type } \alpha \} $.
\end{itemize}
\end{lem}
\begin{demo}
  Statement (1) is immediate from the definition of $ E_A$ given in (\ref{runsoverpairs}), so
  let us prove (2). Suppose that $ A=\{I_1, I_2, \ldots, I_k \} $. Then by construction, the vector
  $  v_{\underline{i}}^{\underline{s}} $ belongs to the right hand side of (2) exactly when any two 
  terms $ s_i$ and $s_j$ of $ \underline{s}$ coincide iff $ i $ and $ j $ are in the same block $ I_l$ of $ A$.
  For example, if $ r =3 $ and
  $ A= \{ \{ 1,2\},  \{ 3,4\},   \{ 5,6,7\} \}$ then the $ \underline{s}$ satisfying
  $ A_{\underline{s}} = A$ are the following ones 
  \begin{equation}
    \begin{array}{l}
(1,1,2,2,3,3,3), \,(2,2,1,1,3,3,3), \,(3,3,1,1,2,2,2), \,(1,1,2,2,3,3,3), \\(1,1,3,3,2,2,2), (2,2,3,3,1,1,1).
\end{array}
  \end{equation}
  Hence, in view of (\ref{moebius}), that is 
\begin{equation}\label{key}
  \mathbb{E}_A= \sum_{A \subseteq B }\mu(A,B)E_B, 
\end{equation}
together with 
  $ \mu(A,A) = 1 $, we conclude from (1) 
  that $ \mathbb{E}_A  v_{\underline{i}}^{\underline{s}} =  v_{\underline{i}}^{\underline{s}} $
  whenever $  v_{\underline{i}}^{\underline{s}} $ belongs to the right hand side of (2), which proves the inclusion
  $ \supseteq $ of (2). To prove $ \subseteq $, we assume that
  $ v \in \be_A   V^{\otimes n} $ and consider its expansion $ v = \sum
  \alpha_{\underline{i}}^{\underline{s}}  v_{\underline{i}}^{\underline{s}}$ with coefficients 
  $ \alpha_{\underline{i}}^{\underline{s}} \in S$. By (1) and (\ref{key})
  we have that $ A \subseteq A_{\underline{s}} $
  whenever $ \alpha_{\underline{i}}^{\underline{s}} \neq 0 $ and hence
  if $ v $ does not belong to the right hand side of (2), there are
  $\underline{i}_0 \in \seqN $ and $ \underline{s}_0 \in \seqr $ 
  such that $ \alpha_{\underline{i}_0}^{\underline{s}_0} \neq 0 $ and
  $ A  \subsetneqq A_{\underline{s}_0} $.
  We have that $ E_{A_{\underline{s}_0}} v_{\underline{i}_0}^{\underline{s}_0}  = v_{\underline{i}_0}^{\underline{s}_0}  $
      and hence $ E_{A_{\underline{s}_0}} v \neq 0 $, by (1) once again. 
   On the other hand, by
   part (3) of Proposition \ref{commuting} we also have that $ E_{A_{\underline{s}_0}} v= 0 $,
   and this gives the desired
   contradiction.
Finally, (3) follows from (2) and the definitions.

\end{demo}

\medskip
For $ N $ any natural number we define 
$ \MCN \subseteq \MC $ (resp. $ \MPN \subseteq \MP $) as the set of multicompositions
(resp. multipartitions) 
$\blambda=(\lambda^{(1)},\lambda^{(2)},\ldots,\lambda^{(r)})$ such
that each $ \lambda^{(i)} $ is of length less than $ N$.

\medskip
Suppose now that 
$ \blambda \in \MCN $ and that $ \bT \in \rstd(\blambda) $. 
Then we define $ \underline{i}^{\bT}= (i_1, \ldots, i_n) \in \seqN$
by letting $ i_j $ be the row number of the node in $ \bT $ containing $ j$,
and similarly we define $ \underline{s}^{\bT}= (s_1, \ldots, s_n) \in \seqr$
by letting $ s_j $ be the component number of the node in $ \bT $ containing $ j$.
For example if
\begin{equation}\label{newaboveexample}
\bT={\footnotesize\left(\;\young(3,4)\,,\,\young(1,8)\,,\,\young(69)\,,\,\young(7\diez)
  \,,\,\ytableausetup{centertableaux,boxsize=1.3em}\begin{ytableau}
    11\\12\\13\end{ytableau}\,,\,\begin{ytableau}14\\15\\16
    \end{ytableau}\,,\,\begin{ytableau}20\\18\\19
    \end{ytableau}\;\,,\,\begin{ytableau}23 &21\\ 24
    \end{ytableau}\;\,,\;\begin{ytableau} 17 \\ 22 & 25
    \end{ytableau}\;, \;\young(5,2)\right)} 
\end{equation}
then $ \underline{i}^{\bT} = (1,2,1,2,1,1,1,2,1,1,1,2,2,1,2,3,1,2,3,1,1,2,1,2,2)$
and $ \underline{s}^{\bT}= (2,10,1,1,10,3,4,\newline2,3,4,5,5,5,6,6,6,9,7,7,7,8,9,8,9)$.
For $ \bT \in \rstd(\blambda)$ we write $ v_{\bT} $ for $v_{\underline{i}^{\bT}}^{\underline{s}^{\bT}}$.

\medskip
The following Theorem relates the tensor space module $ V^{\otimes n} $
with the permutation modules $ M(\blambda) $. 

\begin{teo}\label{17}
Let $ V $ be the free $ S$-module (or $ \kk $-vector space) of dimension $ rN $ with basis as in 
(\ref{95}). Suppose that 
$ \blambda \in \MCN $ is of type $ \alpha \in \Par$ and that $ \alpha $ is of length $ r $.
  \begin{itemize}
\setlength\itemsep{-0.5em}
\item[(1)] The $S$-linear map $\iota_{\blambda}: M(\blambda) \rightarrow V^{\otimes n} $
  given by $\iota_{\blambda}(x_{\bT}) = v_{\bT} $ for $ \bT \in \rstd(\blambda) $, is
  an embedding of $ \Ea$-modules. 
  \newline
\item[(2)] Identifying $ M(\blambda) $ with $ \varphi_{\blambda}( M(\blambda))$ we have
 $ M(\blambda) \subseteq  \be_{\alpha}V^{\otimes n} $ and 
  a direct sum decomposition of $ \Ea $-modules
  \begin{equation}\label{directsun}
   \be_{\alpha} V^{\otimes n} \cong \bigoplus_{\substack{\blambda \in \MCN\\ |A_{\blambda}|  =\alpha}} M(\blambda).
    \end{equation}
\end{itemize}
\end{teo}
\begin{demo}
  (1) follows from Lemma \ref{descr} together with (\ref{operatorE}) and (\ref{operatorG}) and the definitions
  and (2) follows from (1) and Lemma \ref{16}.
\end{demo}

\medskip
For any natural number $ N $ we define $ {\cal L}_{n, \le N}(\alpha) $ as the
subset of $  {\cal L}_n(\alpha) $ consisting of the pairs $ (\blambda \mid \bmu )  $ such
that the components $\lambda^{(k)}$ of $ \blambda $ all have less than $ N $ columns, and
we set $ {\cal L}_{n, \nleq N}(\alpha) :={\cal L}_{n}(\alpha) \setminus
{\cal L}_{n, \le N}(\alpha) $.
We can now give the main Theorem of our paper.
In the Hecke algebra case the argument of the
proof was discovered in \cite{Har},  although our presentation of the argument, even in
the Hecke algebra case, 
differs slightly from the one {\color{black}used} by H\"{a}rterich.
\begin{teo}\label{mainT}
Let $ V $ be as in (\ref{95}) with basis 
$ \cal B $ and suppose that $ \alpha \in \Par $ is of length $ r$.
Let $ {\mathcal I} \subseteq \Ea$ be the annihilator ideal of the action of $ \Ea $ in
$    \be_{\alpha}V^{\otimes n} $.
  Then $ {\mathcal I} $ is free over $ S $ with basis
  \begin{equation}\label{110}
    \left\{ n_{\es\et} \, | \, \es, \et \in \std(\Lambda),
    \Lambda \in {\cal L}_{n, \nleq N}(\alpha)   \right\} .
  \end{equation}
A similar statement holds over $ \kk$.  
\end{teo}
\begin{demo}
Let us focus on the $S$-case, since the $ \kk $-case is done the same way.
  Let $ {\mathcal I}_1 \subseteq \Ea$ be the $S$-span of
  $\left\{ n_{\es\et} \, | \, \es, \et \in \std(\Lambda), \Lambda \in {\cal L}_{n, \nleq N}(\alpha) \right\}$.
  We first observe that $ {\mathcal I}_1 $ is a two-sided ideal in $ \Ea$, as one
  sees from the $  n_{\es\et} $-version of Lemma 56 of \cite{ER}: the number of columns of the components
  never decreases under the straightening procedure of that Lemma.
We first prove that $ {\mathcal I}_1 \subseteq {\mathcal I} $, that is 
for a basis element $ n_{\es\et} $
of $ {\mathcal I}_1 $ we prove that
$  M(\blambda) n_{\es\et} = 0$ for any $ M(\blambda ) $ appearing
in (\ref{directsun}), or equivalently that $  M(\blambda) n_{\es\et} = 0$.
But since   
$ {\mathcal I}_1 $ is an ideal, in order to show $  M(\blambda) n_{\es\et} = 0$
it is enough to show $ x_{\blambda} n_{\es\et} = 0$. 
Note that $ \blambda $ may not be a multipartition, only a multicomposition, 
but in any case $ \bT^{\blambda} $ is of the initial kind and so we
can use Lemma 4.4 of \cite{Mur95} to rewrite $ x_{\blambda} $
as a linear combination
of $ x_{\Bs_1 \bT_1 } $ where $ \Bs_1 $ and $ \bT_1  $ are $\blambda_1 $-multitableaux of the initial kind, 
for $ \blambda_1 $ an $ r$-multipartition such that
$ \blambda_1 \unrhd \blambda $. In general,
for any $ \sigma \in \Si_r$ there is an $ \Ea$-isomorphism
$ M(\blambda) \cong M(\blambda^{\sigma})  $ 
where $ \blambda^{\sigma} $ is the $r$-multicomposition obtained
from $ \blambda^{\sigma} $ by permuting the components,
as one checks using
Lemma \ref{descr}, 
and so we may 
assume that $ \blambda_1 $ is increasing. Moreover, since $ \blambda_1 \unrhd \blambda$ 
and $ \blambda \in \MCN $ we have $ \blambda_1  \in   \MPN   $.
We now observe that 
$ x_{\Bs_1 \bT_1 } $ is also of the form $ m_{\es_1 \et_1} $, where
$ \es_1 = ( \Bs_1 \mid  \bu_1) $ and $ \et_1 = ( \bT_1 \mid  \bv_1) $ and 
where the components of $ Shape(\bu_1) $ and $ Shape(\bv_1) $ are all one-column partitions,
in particular $ d(\bu_1) = d(\bv_1) = 1$.
But from Lemma {\ref{secondduality}} we know that 
$  m_{\es_1 \et_1} n_{\es\et} \neq  0 $ implies
$ \et_1 \unlhd \es^{\prime} $
which is impossible since all
components of $  \blambda_1 $ have less then $ N$ rows, whereas at least one component of 
$ Shape(\Bs^{\prime}) $ has more then $ N $ rows; here $ \es = (\Bs \mid \bu)$ for some $ \bu$.
Hence
$  m_{\es_1 \et_1} n_{\es\et} =  0 $ and so 
the inclusion $ {\mathcal I}_1 \subseteq {\mathcal I} $
is proved.

\medskip
Suppose now that the inclusion $ {\mathcal I}_1 \subseteq \mathcal I$ is strict. Then there
exists $ n \in { \mathcal I} \setminus {\mathcal I}_1 $.
Our plan is  then to construct a $ \blambda  \in \MCN $
such that $ | A_{\blambda} | = \alpha$ and
$ M(\blambda) n \neq 0 $,
which is a contradiction since $ M(\blambda) $ is a direct summand of
$ \be_{\alpha} V^{\otimes n} $ by Lemma \ref{17}. We shall do so by showing that there exists
a $ \Lambda := (\blambda \mid \bmu) \in {\cal L}_{n}(\alpha) $ with  $ \blambda  \in \MCN   $,
such that $ M(\Lambda) n \neq 0 $. This 
gives the desired contradiction since $ M(\Lambda)  \subseteq M(\blambda) $. 

In order to do so we consider the  
expansion $ n = \sum_{\es, \et \in \std(\Lambda), \Lambda \in {\cal L}(\alpha)} \lambda_{\es \et}n_{\es \et} $.
We may assume that $ Shape(\es) \in {\cal L}_{n, \le N}(\alpha) $ for all $ (\es, \et) $ 
occurring in the expansion of $ n $: otherwise if 
$ Shape(\es) \in {\cal L}_{n, \nleq N}(\alpha) $
we subtract the corresponding term $ \lambda_{\es \et}n_{\es \et} $ from $ n$ and
still get an element in $ {\mathcal I}_1 \subseteq \mathcal I$, by the first part of the proof.
We now choose $ (\es_{min}, \et_{min}) $ with 
$ \lambda_{\es_{min} \et_{min}} \neq 0$ and minimal in the sense that 
if $ \es \unlhd \es_{min}$, 
$ \et \unlhd \et_{min}$ and $ \lambda_{\es \et} \neq 0$ 
then $ (\es, \et ) = (\es_{min}, \et_{min} ) $. Let $ \es_{min} = (\Bs_{min } \mid \bu_{min}) $
and let $ \Lambda_{min}= Shape(\es_{min}) = (\blambda_{min}\mid \bmu_{min}) $.
We then have $ \blambda_{min}^{\prime}   \in \MCN $ of type $ \alpha$ 
and so we may consider the permutation module $ M(\Lambda^{\prime}_{min})$.
We are interested in the 
elements $ m_{\es_{min}^{\prime}} $ and $ m_{\et_{min}^{\prime}} $ of $ M(\Lambda^{\prime}_{min})$.
Using the minimality of $ (\es_{min}, \et_{min}) $, together with Lemma
{\ref{secondduality}} and the invariance of the bilinear form, we get
\begin{equation}
  (m_{\es_{min}^{\prime}} n , m_{\et_{min}^{\prime}})_{\Lambda^{\prime}} =
    (m_{\es_{min}^{\prime}} n^{}_{ \es^{}_{min} \et^{}_{min}}    , m_{\et_{min}^{\prime}})_{\Lambda^{\prime} }
\end{equation}
which is nonzero by Lemma \ref{samespirit}. Hence
$ M(\Lambda^{\prime}_{min}) n \neq 0 $, as needed. This proves the Theorem.
\end{demo}

\medskip
In view of the above Theorem, it is natural to consider the following quotient algebra
\begin{equation}
\EaTL: =\Ea/\mathcal I
\end{equation}
where $ \mathcal I $ is the ideal from (\ref{110}). It is a generalization of
the 'generalized Temperley-Lieb' algebras introduced in \cite{Har}.
For $ \kk $ a field, there is also a specialized version of $ \EaTL $ that we denote
$ \EaTLK $.
By construction $ \EaTL $
is the largest quotient of $ \Ea $ 
acting faithfully on
$  \be_{\alpha} V^{\otimes n}$. We have the following Corollary to Theorem \ref{mainT}.

\begin{coro}\label{cormaun}
$ \EaTL $ is free over $ S $ with basis 
$    \left\{ n_{\es\et} \, | \, \es, \et \in \std(\Lambda),
  \Lambda \in  {\cal L}_{n, \le N}(\alpha)  \right\} $. A similar statement holds
  for $ \EaTLK$.
\end{coro}
The following Corollary shows that 
$ \EaTL $ is also a generalization of the \textit{partition Temperley-Lieb algebra} $ \PTLa$ 
that was introduced in \cite{Juyumaya}. 

\begin{coro}
  Let $ St_i \in \Ea$ be the $ i$'th \textit{Steinberg element} given by 
  \begin{equation} St_i:=-q^{-3}g_i g_{i+1} g_i +q^{-2}g_i g_{i+1} +q^{-2}g_{i+1} g_i - q^{-1}g_i - q^{-1}g_{i+1} +1 \end{equation}
  and define $ \PTLa $ as the quotient algebra $ \PTLa := \Ea/{\cal J}_i $
  where $ {\cal J}_i $ is the two-sided ideal of $ \Ea $ generated by $ e_i e_{i+1} St_i  $, for some $ i$.
  Then $  \PTLa $ does not depend on the choice of $ i $ and, moreover,
$ \PTLa   = \EaTLtwo $.
  Similar statements hold over $ \kk$.
\end{coro}
\begin{demo}
  In Proposition 4.5 of \cite{Juyumaya} it is shown that the $ e_i e_{i+1} St_i $'s are all conjugate
  to each other in $ \Ea $ and so the $ {\cal J}_i $'s are equal ideals, which proves the first statement. 
For simplicity we write $ {\cal J}= {\cal J}_1   $. 
  We must then show that $ {\cal J }=  {\cal I } $ where $ \cal I$ is the ideal
  described in Theorem \ref{mainT}. The inclusion $ {\cal J } \subseteq  {\cal I } $ follows from
  the fact that $ e_{n-1} e_n St_{n-1} $ belongs the basis for $ \cal I $ given in (\ref{110}). In fact we have
  $ e_{n-1} e_n St_{n-1}  = n_{\et^{\Lambda} \et^{\Lambda}}$ where $ \Lambda = 
  ( \blambda \mid \bmu ) $ is chosen such that all components of $ \blambda = 
  ( \lambda^{(1)}, \ldots , \lambda^{(r)}) $ are one-column partitions, except $ \lambda^{(r)} $
  which is of the form $ \lambda^{(r)} = (3,1^s) $ for some $ s$, and such that all
  components of $ \bmu $ are one-column partitions. For the other inclusion
  $ {\cal J } \supseteq  {\cal I } $ we consider a basis element $ n_{\es\et}  $ for $ \cal I $, as given in
  (\ref{110}). Letting $ Shape(\es) = (\blambda \mid \bmu ) $, 
  there is a component of $ \blambda $ with more than three columns and so there is an
  $ i$ such that 
  $ e_{i} e_{i+1} St_{i} $ is a factor of $ n_{\es\et}  $. This shows
  $ {\cal J } \supseteq  {\cal I } $ and 
  concludes the proof of the Corollary.
\end{demo}

\medskip
\noindent
    {\bf Remark}. Using Corollary \ref{cormaun} we can determine the
    dimension of $ \EaTL $, and hence also of $ \PTL:= \bigoplus_{\alpha \in \Par} \PTLa$,
    since it is the cardinality of
    \begin{equation}\label{correct} \bigcup_{\alpha \in \Par} \left\{ n_{\es\et} \, | \, \es, \et \in \std(\Lambda),
        \Lambda \in  {\cal L}_{n, \le 2}(\alpha)  \right\}.
    \end{equation}
        We have for example
        \begin{equation}
\dim \PTLthree =29 , \, \, \,\, \, \,\, \, \,\, \, \, 
          \dim \PTLfour =334 , \, \, \,\, \, \,\, \, \,\, \, \,    \dim \PTLfive =    5512.  \end{equation}
        These values are confirmed by Juyuamaya and Papi's MAGMA calculations and
        by Espinoza's preprint, \cite{Esp}, that contains a closed general formula for $ \dim \PTL$.
    
\medskip
\noindent
    {\bf Remark}. In Conjecture 6.5 of \cite{Juyumaya} there is a conjectural basis for $  \PTL$
    which would imply that $ \dim  \PTL = b_n c_n$ where $ c_n = \frac{1}{n+1} {2n \choose n}$
    is the Catalan number. It follows from (\ref{correct}) that Juyumaya's conjectural basis is wrong, as his and
    Papi's calculations had in fact already shown.  

\medskip
\noindent
    {\bf Remark}.
It would be interesting to use (\ref{correct}) to
    study Conjecture 7.3 from \cite{Juyumaya} on the existence of a Markov trace on $ \PTL $.
    One would expect a diagrammatic calculus associated with $\PTL$ that might be of relevance for this.

\noindent
\sc Instituto de Matem\'aticas, Universidad de Talca, Chile,
steen@inst-mat.utalca.cl, steenrh@gmail.com

\end{document}